\documentclass{elsart} 
\usepackage{latexsym}
\usepackage{amsfonts}
\usepackage{amssymb}
\usepackage{amsmath}
\usepackage{amscd}
\usepackage{eucal}
\usepackage[all]{xy}
\usepackage{pstricks}
\newcommand{\bs}{\boldsymbol}
\newcommand{\mb}{\mathbf}

\begin{document}
\theoremstyle{definition} \newtheorem{Defi}{Definition}
\theoremstyle{definition} \newtheorem{nota}{Notation}



\begin{frontmatter}
\title{Loop spaces, and coherence for monoidal and braided monoidal bicategories}
\author[ng]{Nick Gurski\thanksref{me}}
\address[ng]{Deparment of Pure Mathematics, University of Sheffield}
\ead{nick.gurski@sheffield.ac.uk}
\thanks[me]{The author would like to thank Scott Carter, Andrew Casson, Eugenia Cheng, and Peter May for enlightening conversations.}



\begin{abstract}
We prove a coherence theorem for braided monoidal bicategories and relate it to the coherence theorem for monoidal bicategories.  We show how coherence for these structures can be interpretted topologically using up-to-homotopy operad actions and the algebraic classification of surface braids.
\end{abstract}
\end{frontmatter}


\section{Introduction}

Braided monoidal categories have appeared in various branches of mathematics in the past two decades.  They show up in category theory as the centers of monoidal categories, and in higher dimensional category theory as special kinds of weak 3-categories \cite{gps}.  Categories whose morphisms are of a geometric nature are often braided, such as categories of ribbons.  Many categories related to representation theory often have nontrivial braiding.  For instance, complex representations of $GL_{n}(\mathbb{F}_{q})$, where here $n$ is allowed to vary, give rise to a graded ring by taking characters.  This ring is commutative, but the tensor product of representations from which the multiplication is derived is merely braided \cite{js2}.

Braided monoidal 2-categories are a much more recent subject of study, and sadly fewer concrete examples have been constructed.  The first definition was given by Kapranov and Voevodsky in \cite{kv2}.  This definition was later improved upon by Baez and Neuchl in \cite{bn} by requiring that the two solutions to the Zamalodchikov tetrahedron equation, $S_{+}$ and $S_{-}$, be equal.  Sjoerd Crans later made a further improvement by adding new unit conditions to the definition in \cite{cr}.

It is important, though, that braided monoidal 2-categories were defined and studied before braided monoidal bicategories.  This is likely due to computational issues, most notably that performing computations in a monoidal 2-category is far simpler than performing them in a monoidal bicategory; adding a braiding only serves to make working in the semi-strict setting even more attractive.  Monoidal 2-categories are honest monoids in the monoidal category $\mb{Gray}$ of 2-categories equipped with the Gray tensor product, while monoidal bicategories are essentially one-object tricategories and thus only monoids in a very weak sense.  Monoidal bicategories are much more difficult to deal with computationally than $\mb{Gray}$-monoids, although the coherence theorem in \cite{g1} goes a long way towards remedying this difference.

Every monoidal bicategory is monoidally biequivalent to a monoidal 2-category (that is, a $\mb{Gray}$-monoid) by the coherence theorem for tricategories \cite{gps,g1}, but there are deeper questions about \textit{braided} monoidal bicategories which have no analogue in the realm of braided monoidal categories.  The key insight here is that the only strictness one can impose on a braided monoidal category is that the underlying monoidal category be strict, the braiding itself can neither be strict nor non-strict.  There is a definition of a symmetric monoidal category which does impose an extra axiom, but this gives a new structure, not an equivalent but more well-behaved version of a braided monoidal category.  The difference at the two-dimensional level is that the definition of a braided monoidal 2-category asks that the braiding be somewhat strict by imposing axioms on some of the 1-cells involved.

Thus we arrive at two crucial coherence questions:  what is the weakest possible definition of a braided monoidal bicategory, and how does this structure relate to the braided monoidal 2-categories of previous authors?  The answer to the first question turns out to be relatively straightforward.  Weakening the definition of Baez and Neuchl is a simple matter of altering pasting diagrams, as is weakening the definition of Crans.  The weakening \textit{procedure} is largely uninteresting by the coherence theorem for monoidal bicategories, but the final product is far less trivial.  For instance, our weakening of the Baez-Neuchl definition completely ignores Crans' unit axioms, but the coherence results we prove -- Theorems \ref{cohbmb} and \ref{strictbmb} in particular -- recover them in the strictification process.

The central theoretical concern of this paper is understanding the structure of free braided monoidal bicategories.  Just as in studying the coherence theory for monoidal and braided monoidal categories, understanding what equations hold in free braided monoidal bicategories is the first step to being able to easily calculate within them.  We solve this problem completely in the case of the free braided monoidal bicategory on a set $S$ by showing that in this case a pair of parallel 1-cells is either uniquely isomorphic or not isomorphic at all.

This coherence theorem highlights a key feature of both braided monoidal categories and braided monoidal bicategories:  the interesting structure occurs on the level of 1-cells.  In the braided monoidal category case, this is the top dimension of the structure.  But in the braided monoidal bicategory case, we also have 2-cells, and the coherence theorem shows that these 2-cells do not introduce any new braided phenomena that was not already present at the 1-cell level.  This is the reflection in category theory that the configuration space of $n$ unordered points in the plane is a $K(B_{n},1)$ where $B_{n}$ is the $n$th braid group.  Thus while there are many non-isomorphic parallel 1-cells in this free object corresponding to different elements of the fundamental group of the configuration space, if two parallel 1-cells are isomorphic then they are uniquely so, corresponding to the fact that $\pi_{2}$ of this configuration space is zero.

While the categorical and bicategorical theories produce very similar structures, there is a significant increase in the complexity of the algebra by moving up a dimension:  the free braided monoidal category on one object is small enough to understand directly in a generators-and-relations fashion, while the same cannot be said of the corresponding free braided monoidal bicategory.  It is relatively simple to show that the free braided monoidal category on one object is equivalent, as a braided monoidal category, to the disjoint union of the braid groups now thought of as one-object groupoids.  This disjoint union $\mathcal{B}$ has objects the natural numbers $[n]$ and hom-sets empty if the source and target are different or equal to the $n$th braid group in the case of $\mathcal{B}([n], [n])$.  The braided monoidal structure on $\mathcal{B}$ is quite easy to describe directly by drawing pictures of braids.  This equivalence of categories makes understanding free braided monoidal categories much simpler, and the resulting coherence theorems allow computations in an arbitrary braided monoidal category to be greatly simplified.  While a similar result is true for braided monoidal bicategories, the algebra is much too difficult to tackle directly.

To understand the algebra of the free braided monoidal bicategory on one object, it is instructive to first understand the algebra of the free braided monoidal category on one object in a more sophisticated way.  First, recall that the $n$th braid group can be defined as the fundamental group of the configuration space of $n$ unordered points in $\mathbb{R}^{2}$, denoted $B(n, \mathbb{R}^{2})$.  This space is path-connected, so we can replace the fundamental group by the fundamental groupoid as these are equivalent as categories.  Now the fundamental groupoid functor commutes with coproducts, and thus the free braided monoidal category on one object is equivalent to
\[
\Pi_{1} \Big( \coprod B(n, \mathbb{R}^{2}) \Big).
\]

The braided monoidal structure on this category is still somewhat artificial and to remedy this we utilize the theory of symmetric operads.  First we note that $ B(n, \mathbb{R}^{2})$ is homotopy equivalent to $\mathcal{C}_{2}(n)/\Sigma_{n}$, where $\mathcal{C}_{2}$ is the little 2-cubes operad, so
\[
\coprod B(n, \mathbb{R}^{2}) \simeq \coprod \mathcal{C}_{2}(n)/\Sigma_{n}.
\]
But recall that $\coprod \mathcal{C}_{2}(n)/\Sigma_{n}$ is actually the free $\mathcal{C}_{2}$-algebra on a terminal space.  It is easy to show that if a space $X$ is a $\mathcal{C}_{2}$-algebra, then $\Pi_{1}X$ is actually a braided monoidal category.  Using the homotopy equivalence above, it is then straightforward to show that $\mathcal{B}$ is equivalent, as a braided monoidal category, to the  braided monoidal category $\Pi_{1}\Big( \coprod \mathcal{C}_{2}(n)/\Sigma_{n} \Big)$ obtained using the algebra structure over $\mathcal{C}_{2}$.  Therefore we have replaced the very hands-on definition of the braided monoidal structure on $\mathcal{B}$ with the topologically-derived braided monoidal structure on $\Pi_{1} \Big( \coprod \mathcal{C}_{2}(n)/\Sigma_{n} \Big)$.

The coherence theorem for braided monoidal bicategories follows the same strategy.  We show that the free braided monoidal bicategory on one object is biequivalent, as a braided monoidal bicategory, to the braided monoidal bicategory
\[
\Pi_{2} \Big(  \coprod \mathcal{C}_{2}(n)/\Sigma_{n} \Big),
\]
where the braided monoidal structure on this fundamental 2-groupoid is understood completely in terms of the $\mathcal{C}_{2}$-algebra structure.  Analyzing which equations of 2-cells hold in the free object amounts to understanding the topology of configuration spaces and the little 2-cubes operad.

The topology of these configuration spaces is intimately linked with the geometry of surfaces embedded in $\mathbb{R}^{4}$.  The 2-cells of the bicategory
\[
\Pi_{2} \Big( \coprod B(n, \mathbb{R}^{2}) \Big)
\]
are homotopy classes of homotopies between braids, which after smooth approximation can be represented by smooth surfaces embedded in $\mathbb{R}^{4}$.  These surfaces have boundary, and the pieces of the boundary coming from the source and target of these 2-cells are both braids in $\mathbb{R}^{3}$.  Carter and Saito \cite{cs2} have studied a very similar classification problem, and have classified ambient isotopy classes of certain surfaces embedded in $\mathbb{R}^{4}$ under the name of \textit{braid movie moves}.  The surfaces themselves are the braid movies, and the braid movie moves give an algebraic description of possible ambient isotopies between them.  Because our surfaces are homotopies and not general braid movies, it turns out that the relevant braid movie moves are precisely what they call the C-I braid movie moves.  Our homotopy classes correspond to Carter and Saito's ambient isotopy classes, so their classification theorem then becomes a complete description of the equations between 2-cells in this fundamental 2-groupoid.  But the C-I braid movie moves of Carter and Saito also give the axioms for braided monoidal bicategories (modulo strictness questions which are handled by coherence for monoidal bicategories), thus establishing a braided monoidal biequivalence between the free braided monoidal bicategory on one object and this fundamental 2-groupoid.

The strategy employed here has obvious generalizations in two different directions corresponding to the categorical dimension and the topological codimension.  Studying $\Pi_{n} \Big( \coprod B(k, \mathbb{R}^{2}) \Big)$ should yield information about braided monoidal $n$-categories, although this is far beyond the scope of current technology in higher dimensional category theory.  It is also possible to increase the topological codimension by changing $\mathbb{R}^{2}$ to $\mathbb{R}^{3}$, $\mathbb{R}^{4}$, or higher values.  One might hypothesize that the case of $\mathbb{R}^{3}$ should be related to \textit{sylleptic} monoidal bicategories in a fashion analogous to the braided case.  Here we would use the little 3-cubes operad to induce a sylleptic structure on the fundamental 2-groupoid.  The roadblock now is geometric rather than categorical:  there is no classification theorem for surface braids in $\mathbb{R}^{5}$ whose source and target braids lie in $\mathbb{R}^{4}$.  The extra structure giving a sylleptic monoidal bicategory is an isomorphism $\gamma^{2} \cong 1$ between the braiding squared and the identity.  This isomorphism could certainly be realized geometrically, but the Morse-theoretic classification of ambient isotopies between these types of surface braids has not yet been carried out.

The paper proceeds as follow.  Section 2 is an overview of relevant background material, including braided monoidal bicategories and 2-categories, configuration spaces, surface braids, and symmetric operads.  This section will not have detailed proofs.  For a discussion of configuration spaces and braids, we refer the reader to \cite{fan} and \cite{fon}, and for a treatment of symmetric operads in the topological case the original work of May \cite{m} is a good reference.  The reader should consult the books of Kamada \cite{kam} or Carter and Saito \cite{cs2} for more details of the results needed about surface braids.

Section 3 gives the full construction of the fundamental 2-groupoid.  We first construct a tricategory of topological spaces, and then show that the fundamental 2-groupoid is a functor of tricategories $\Pi_{2}:\mb{Top} \rightarrow \mb{Bicat}$.  We then show that when the spaces involved are algebras over the little $1$- or $2$-cubes operads, then the fundamental 2-groupoids can be equipped with the structure of a monoidal or braided monoidal bicategory, respectively.  These structured fundamental 2-groupoids then allow us to equip some bicategories with extra structure by using the fact that the forgetful functors
\[
\begin{array}{c}
\mb{MonBicat} \rightarrow \mb{Bicat},\\
\mb{BrMonBicat} \rightarrow \mb{Bicat}
\end{array}
\]
both lift biadjoint biequivalences.  Thus even spaces which are only homotopy equivalent to algebras for the operads $\mathcal{C}_{n}$ have fundamental 2-groupoids which come equipped with additional monoidal structure once a homotopy equivalence is chosen.  In particular, the fundamental 2-groupoid of the space $\coprod B(n, \mathbb{R}^{2})$ is braided monoidal, even though this space is not itself an algebra for $\mathcal{C}_{2}$.

Section 4 studies free monoidal and braided monoidal bicategories in detail.  This is where we prove our coherence theorems.  This section splits into two parts.  For the first, we use the already known coherence theory for monoidal bicategories to show that the free such object is appropriately biequivalent to the fundamental 2-groupoid of a coproduct of configuration spaces.  In the second we reverse this, proving the biequivalence and then using that to give a coherence result using topology and geometry.

Section 5 gives a strictification result.  We show that every braided monoidal bicategory is appropriately biequivalent to a braided monoidal 2-category in the sense of Crans \cite{cr}.  Here we assume the reader is familiar with the strictification construction given for tricategories in \cite{g1}.

We note here that all concepts are the maximally weak version unless stated otherwise, so that functor means weak functor or weak monoidal functor, etc.  Thus the prefix ``2-'' will always be used to indicate strictness with one exception: we use the phrase ``fundamental 2-groupoid'' even though the construction is weak.  The reader should be cautioned, though, that monoidal 2-categories are not monoids in the category $\mb{2Cat}$ with the Cartesian product, but rather they are monoids in $\mb{Gray}$ which is the same category but with the monoidal structure given by the Gray tensor product.

\section{Background}

This section will give the necessary background definitions and results for the rest of the paper.  We begin by reviewing the definition of the little $n$-cubes operad $\mathcal{C}_{n}$ and stating some basic results.  Every operad gives rise to a monad, and we write the monad associated to $\mathcal{C}_{n}$ as $C_{n}$.  Then we recall how certain free algebras over $C_{n}$ are $\Sigma_{n}$-equivariantly homotopy equivalent to configuration spaces.  Next we review the definitions and basic results of surface braid theory.  Finally, we give the central definitions of braided monoidal bicategory and braided monoidal 2-category (following Crans \cite{cr}).

\subsection{Operads and monads}

Recall that the data for a symmetric operad $\mathcal{P}$ in a symmetric monoidal category $M$ consists of objects $\mathcal{P}(n)$ for all integers $n \geq 0$, maps
\[
\alpha_{n, (k_{1}, \ldots, k_{n})}: \mathcal{P}(n) \otimes \mathcal{P}(k_{1}) \otimes \cdots \otimes \mathcal{P}(k_{n}) \rightarrow \mathcal{P}(k_{1} + \cdots + k_{n}),
\]
a unit map $\iota:I \rightarrow \mathcal{P}(1)$, and a right $\Sigma_{n}$-action on $\mathcal{P}(n)$ for each $n$.  These maps must satisfy associativity, unit, and equivariance axioms.  These axioms are actually just the axioms for the sequence of objects $\{ \mathcal{P}(n) \}$ to be a monoid in the monoidal category of symmetric collections in $M$.

If $\mathcal{P}$ is a symmetric operad in $M$,then an algebra for $\mathcal{P}$ consists of an object $X \in M$ and maps
\[
\mathcal{P}(n) \otimes X^{\otimes n} \rightarrow X
\]
satisfying associativity, unit, and equivariance axioms.  If $M$ happened to be closed monoidal, then the symmetric collection defined by $End_{X}(n) = M(X^{\otimes n}, X)$ is an operad, with multiplication given using both composition and the tensor product.  In this case, an algebra structure on $X$ is nothing more than a map of operads $\mathcal{P} \rightarrow End_{X}$.  There is then a category of algebras for $\mathcal{P}$ with morphisms defined to be those maps of underlying objects which strictly preserve the action of $\mathcal{P}$.  We refer the reader to \cite{m} for the complete definitions in the case $M = Top$ as this is the only case needed for our work.

We also have the related notion of a monad on a category.  The data for a monad $(T, \mu, \eta)$ on a category $C$ consists of a functor $T:C \rightarrow C$ and natural transformations $\eta: 1_{C} \Rightarrow T$, $\mu: T \circ T \Rightarrow T$ satisfying associativity and unit axioms.  Once again, these axioms amount to nothing more than the axioms for a monoid in the monoidal category of endofunctors on $C$.

Given a monad $(T, \mu, \eta)$ on $C$, an algebra for $T$ consists of an object $X \in C$ and a morphism $\alpha:TX \rightarrow X$ satisfying associativity and unit axioms.  Just as with the case of an operad, a monad $T$ gives rise to a category of algebras with morphisms $(X, \alpha) \rightarrow (Y, \beta)$ being those morphisms $f:X \rightarrow Y$ in $C$ that strictly preserve the action of $T$.  We refer the reader to \cite{mac} for a full treatment of monads and their algebras.

The basic result relating operads and monads is the following.

\begin{prop}
Let $\mathcal{P}$ be an operad in $M$.  Then the functor $P$ defined on objects by
\[
X \mapsto \coprod_{n \geq 0} \mathcal{P}(n) \otimes_{\Sigma_{n}} X^{\otimes n}
\]
is a monad.  The category of algebras for the operad $\mathcal{P}$ is equivalent to the category of algebras for the monad $P$.
\end{prop}

The operads in this work will naturally act on based spaces, so the description above of the monad associated to an operad needs to be altered.  We will consider an operad $\mathcal{P}$ in the category of unbased spaces, and algebras for it in the category of based spaces.  These algebras will be the algebras for the monad $P$ defined on objects by
\[
X \mapsto \coprod_{n \geq 0} \mathcal{P}(n) \times_{\Sigma_{n}} X^{\times n}/ \sim,
\]
where $\sim$ generates an equivalence relation making basepoint identifications (see \cite{m} for the precise relation).  In the case that we require our operad $\mathcal{P}$ to act on an unbased space, we first attach a disjoint basepoint.

The most important operads for this work are the little $n$-cubes operads $\mathcal{C}_{n}$ introduced by May \cite{m}.  Let $J$ denote the open unit interval.  A \textit{little n-cube} is a linear embedding $\alpha:J^{n} \rightarrow J^{n}$ which is of the form $\alpha = \alpha_{1} \times \alpha_{2} \times \cdots \times \alpha_{n}$ where each $\alpha_{i}$ is a linear map
\[
\alpha_{i}(t) = (y_{i} - x_{i})t + x_{i}, \quad 0 \leq x_{i} < y_{i} \leq 1.
\]
The space $\mathcal{C}_{n}(k)$ is the subspace of $\textrm{Map} \big( (J^{n})^{k}, J^{n} \big)$ consisting of those $k$-tuples of little $n$-cubes which are pairwise disjoint.  (The reader should note that $\mathcal{C}_{n}(0) = *$ as there is a unique empty collection of little $n$-cubes; this will be important later to produce units for our monoidal bicategories.)  There is an obvious operadic multiplication on the spaces of little $n$-cubes given by composition of maps.  It is simple to check that this is a symmetric operad in the category of spaces.  The fundamental result about $\mathcal{C}_{n}$ is the following \cite{m}.

\begin{thm}
1.  Let $X$ be a pointed space.  Then the $n$-fold loop space of $X$, $\Omega^{n}X$, is an algebra for $\mathcal{C}_{n}$.  \\
2.  Let $X$ be a path-connected pointed space.  If $X$ is also an algebra for $\mathcal{C}_{n}$, then $X$ is weakly equivalent to $\Omega^{n}Y$ for some pointed space $Y$.
\end{thm}

\subsection{Configuration spaces}

Given a space $X$ and a natural number $k$, let $\textrm{Config}(k,X)$ denote the space of $k$-tuples of points $(x_{1}, \ldots, x_{k})$ in the space $X$ such that $x_{i} \neq x_{j}$ if $i \neq j$ with the topology induced by the obvious inclusion
\[
\textrm{Config}(k,X) \hookrightarrow X^{k}.
\]
Now $X^{k}$ has a free action of the symmetric group $\Sigma_{k}$, and $\textrm{Config}(k,X)$ inherits this action.

\textbf{Definition}  $B(k,X)$ is the quotient space $\textrm{Config}(k,X)/\Sigma_{k}$.

We refer the reader to \cite{m} for a proof of the following proposition.

\begin{prop}
The free $C_{n}$-algebra $C_{n}(*)$ on the terminal (unbased) space * is homotopy equivalent to
\[
\coprod_{k \geq 0} B(k, \mathbb{R}^{n}).
\]
\end{prop}

There is also a version of this result replacing the terminal set $*$ with any set $S$.  To state the result, we must first define the configuration space of $k$ unordered points in $X$ with labels in the set $S$.  Let $\pi:(X \times S)^{k} \rightarrow X^{k}$ be the map which projects onto the $X$-coordinates.  Consider the subspace $\pi^{-1} \Big( \textrm{Config}(k,X) \Big)$ in $(X \times S)^{k}$.  This is the space of $k$-tuples $\Big( (x_{1}, s_{1}), \ldots, (x_{k}, s_{k}) \Big)$ such that the points $x_{i}$ are all distinct; in particular, two different points can have the same label.  This space has a free action of the symmetric group $\Sigma_{k}$.

\textbf{Definition}  $B(k,X; S)$ is the quotient space
\[
\pi^{-1} \Big( \textrm{Config}(k,X) \Big) /\Sigma_{k}.
\]

We record the next two propositions for future use, they are both simple to prove.

\begin{prop}
The free $C_{n}$-algebra $C_{n}(S)$ on the set $S$ viewed as an unbased, discrete space is homotopy equivalent to
\[
\coprod_{k \geq 0} B(k, \mathbb{R}^{n}; S).
\]
\end{prop}

\begin{prop}
The natural map $B(k,X;S) \rightarrow B(k,X)$ which forgets the labels is a fibration with fiber $S^{k}$.
\end{prop}

\subsection{Braid movie moves and surface braids}

In this section we will introduce the classification of surface braids via braid movie moves of Carter and Saito \cite{cs2}.  This classification result serves to relate the topology of configuration spaces to the algebra of braided monoidal bicategories.  We begin with the basic definitions (see \cite{kam}, Chapter 14).

\textbf{Definition}  Let $D_{i}^{2}$, $i=1,2$, denote a pair of 2-disks, and $pr_{i}:D_{1}^{2} \times D_{2}^{2} \rightarrow D_{i}^{2}$ the projection map to the $i$th disk.  Let $Q_{m}$ be a collection of $m$ interior points of $D_{1}^{2}$.  A \textit{surface braid} $S$ of degree $m$ is an oriented 2-manifold embedded properly and locally flatly in $D_{1}^{2} \times D_{2}^{2}$ such that
\begin{itemize}
\item the restriction of $pr_{2}$ to $S$, $pr_{2}|_{S}:S \rightarrow D_{2}^{2}$, is a branched covering map of degree $m$, and
\item the boundary $\partial S$ is $Q_{m} \times \partial D_{2}^{2}$.
\end{itemize}

\textbf{Definition}  Let $S,S'$ be two surface braids.  Then $S$ is \textit{equivalent} to $S'$ if there is an ambient isotopy $\{ h_{u} \}_{u \in [0,1]}$ satisfying the following conditions:
\begin{itemize}
\item $h_{1}(S) = S'$;
\item for each $u \in [0,1]$, $h_{u}$ is fiber-preserving in the sense that there is a homeomorphism $H_{u}:D_{2}^{2} \rightarrow D_{2}^{2}$ such that $pr_{2} \circ h_{u} = H_{u} \circ pr_{2};$ \textrm{and}

\item $h_{u}$ restricted to $D_{1}^{2} \times \partial D_{2}^{2}$ is the identity for each $u \in [0,1]$.
\end{itemize}
Note in particular that equivalent surface braids necessarily have the same boundary.

There are two simpler classes of surface braids that we will be interested in later.  We define them now.

\textbf{Definition}  Let $S$ be a surface braid of degree $m$.  $S$ is \textit{trivial} if $S$ is equivalent to $Q_{m} \times D_{2}^{2}$, and $S$ is \textit{simple} if the branched covering is simple, i.e., if for every branch point $y$ in $D_{2}^{2}$ there exists a unique singular point $x$ with $pr_{2}(x)=y$, and this singular point has degree 2.

Now that we have defined simple surface braids, we can view them as maps into an extended configuration space as follows.  We define the space $\textrm{Config}_{e}(X,k)$ to be the subspace of $X^{k}$ consisting of those points $(x_{1}, \ldots, x_{k})$ such that there exists at most one pair of indices $s < t$ for which $x_{s} = x_{t}$.  The symmetric group acts on this space, so we define the extended configuration space $B_{e}(X,k)$ to be the quotient $\textrm{Config}_{e}(X,k)/\Sigma_{k}$.  It should be clear that the usual configuration space $B(X,k)$ is the subspace consisting of those points for which there is no pair of indices $s < t$ with $x_{s} = x_{t}$; the complement of $B(X,k)$ is called the singular locus, and is denoted $\Sigma_{m}^{(1)}(X)$ by Kamada.  In the case that $X$ is the interior of the 2-disk, simple surface braids can be identified with certain kinds of maps $D^{2} \rightarrow B_{e}(X,k)$.  We record the following crucial lemma for later use and refer the reader to \cite{kam} for more discussion of this viewpoint.

\begin{lem}\label{transverse}
Let $S,S'$ be simple surface braids represented by maps $s, s': D^{2} \rightarrow B_{e}(X,k)$.  If $s,s'$ are homotopic via a map $K:I \times D^{2} \rightarrow B_{e}(X,k)$ such that
\begin{itemize}
\item each map $K_{t}$ has the property that $K_{t}(\partial D^{2}) = Q_{m}$, and
\item each map $K_{t}$ intersects $\Sigma_{m}^{(1)}(\textrm{int }D^{2})$ transversely,
\end{itemize}
then $S$ and $S'$ are equivalent surface braids.
\end{lem}

Turning from the geometric approach to surface braids to an algebraic classification, we introduce the theory leading to Carter and Saito's results on braid movie moves \cite{cs2}.

\textbf{Definition}  A \textit{braid movie} is a sequence
\[
(1 = w_{0}, w_{1}, \ldots, w_{k-1}, w_{k}=1)
\]
where each $w_{i}$ is an element in the free monoid generated by symbols $s_{i}^{\pm 1}$, $i=1, \ldots, n-1$, (i.e., a word in these symbols) that satisfies the following condition.  For index $i$, we have that $w_{i} = w_{i-1}$ or that $w_{i}$ differs from $w_{i-1}$ by one of the following elementary braid changes, where here any exponent $\epsilon$ is either $1$ or $-1$:
\begin{enumerate}
\item insertion or deletion of $s_{j}^{\epsilon}$,
\item insertion or deletion of a pair $s_{j}^{\epsilon}s_{j}^{-\epsilon}$,
\item replacement of $s_{j}^{\epsilon_{1}} s_{l}^{\epsilon_{2}}$ with $s_{l}^{\epsilon_{2}} s_{j}^{\epsilon_{1}}$ if $|j - l| > 1$, and
\item replacement of $s_{l}^{\epsilon} s_{j}^{\epsilon} s_{l}^{\epsilon}$ with $s_{j}^{\epsilon} s_{l}^{\epsilon} s_{j}^{\epsilon}$ or replacement of $s_{l}^{\epsilon} s_{j}^{\epsilon} s_{l}^{-\epsilon}$ with $s_{j}^{-\epsilon} s_{l}^{\epsilon} s_{j}^{\epsilon}$ if $|j - l| = 1$.
\end{enumerate}

Every braid movie gives rise to a surface braid by interpreting the elementary braid changes as embedded surfaces in $D_{1}^{2} \times D_{2}^{2}$, and every surface braid can be decomposed into a braid movie.  We now seek to understand how the braid movie representation can help determine when two surface braids are equivalent.  To this end, Carter and Saito \cite{cs1} defined fourteen braid movie moves.  These fall into four groups called the C-I, C-II, C-III, and C-IV moves which are related to the C-moves of Kamada.  We will not define the moves here, but we do note in passing that it is precisely the C-I moves that do not involve branch points.

Before we can state the main result of this section, we must explain one more way to alter a braid movie.  Suppose that $(w_{0}, \ldots, w_{k})$ and $(w_{0}', \ldots, w_{k'})$ are two braid movies such that $w_{j} = w_{j}'$ for all indices except a single one $i$.  Suppose in addition we have that $w_{p} = u_{p}v_{p}$ and $w_{p}' = u_{p}'v_{p}'$ for $p \in \{ i-1, i, i+1 \}$ such that
\begin{itemize}
\item $u_{i}$ is obtained from $u_{i-1}$ by the elementary braid change $\eta$ and
 $v_{i} = v_{i-1}$,
\item $u_{i}' = u_{i-1}'$ and
$v_{i}'$ is obtained from $v_{i-1}'$ by the elementary braid change $\xi$,
\item $u_{i} = u_{i+1}$ and $v_{i+1}$ is obtained from $v_{i}$ by $\xi$, and
\item $u_{i+1}'$ is obtained from $u_{i}'$ by $\eta$ and $v_{i+1}' = v_{i}$.
\end{itemize}
In this case, we say that $(w_{0}', \ldots, w_{k'})$ is obtained from $(w_{0}, \ldots, w_{k})$ by a locality change.

\begin{thm}[Carter and Saito, \cite{cs1}]
Two braid movies represent equivalent surface braids if and only if they are related by \begin{enumerate}
\item a sequence of C-I, C-II, C-III, and C-IV braid movies moves,
\item alterations of braid movie moves via replacing $s_{j}$ with $s_{j}^{-1}$, running any sequence of braid movie move backwards, or replacing a sequence of braid movie moves by the sequence of its palindromes, and
\item locality changes.
\end{enumerate}
\end{thm}

\subsection{Braided monoidal bicategories and 2-categories}

Here we present the definitions of braided monoidal bicategory and braided monoidal 2-category.  A fully weak definition of braided monoidal bicategory has not been proposed in the literature, although the one given here is an obvious weakening of the various definitions of braided monoidal 2-category.  This definition follows the standard philosophy of categorification in which we replace the old axioms (the two hexagons in the definition of a braided monoidal category) with new isomorphisms (the modifications $R_{(-,-|-)}$ and $R_{(-|-,-)}$) and then add new axioms between these.  We also give the definition of braided monoidal 2-category as finalized by Crans in \cite{cr} as a strictified version of our definition of braided monoidal bicategory.

Before giving the definition, we make four notational comments.  First, we have written the tensor product as concatenation to save space.  Second, adjoint equivalences are written $\mathbf{f}$ and have left adjoint $f$, right adjoint $f^{\centerdot}$, invertible unit $\eta_{f}:I \Rightarrow f^{\centerdot}f$, and invertible counit $\varepsilon_{f}:ff^{\centerdot} \Rightarrow I$.  Third, we freely replace 2-cells by their mates under adjoint equivalences without altering the name of the 2-cell  (see \cite{ks} for a discussion of mates and their various properties).  Finally, we have presented the invertible modifications in the definition below by giving their components on objects instead of displaying the source and target transformations explicitly.  The interested reader should find this alternate description easy to construct.

\textbf{Definition} Let $B = (B, \otimes, I, \mathbf{a}, \mathbf{l}, \mathbf{r}, \pi, \mu, \rho, \lambda)$ be a monoidal bicategory.  Then a \textit{braiding} for $B$ consists of
\begin{itemize}
\item an adjoint equivalence $\bs{R}: \otimes \Rightarrow \otimes \circ \tau$ in $\mathbf{Bicat}(B \times B, B)$, where we define $\tau:B \times B \rightarrow B \times B$ to interchange the coordinates;
\item an invertible modification $R_{(-|-,-)}$ as displayed below;
\[
\xy
{\ar^{R1} (0,0)*+{(AB)C}; (15,10)*+{(BA)C} };
{\ar^{a} (15,10)*+{(BA)C}; (45,10)*+{B(AC)} };
{\ar^{1R} (45,10)*+{B(AC)}; (60,0)*+{B(CA)} };
{\ar_{a} (0,0)*+{(AB)C}; (15,-10)*+{A(BC)} };
{\ar_{R} (15,-10)*+{A(BC)}; (45,-10)*+{(BC)A} };
{\ar_{a} (45,-10)*+{(BC)A}; (60,0)*+{B(CA)} };
{\ar@{=>}^{R_{(A|B,C)}} (29, 7)*{}; (29,-7)*{} }
\endxy
\]
\item and an invertible modification $R_{(-,-|-)}$ as displayed below;
\[
\xy
{\ar^{1R} (0,0)*+{A(BC)}; (15,10)*+{A(CB)} };
{\ar^{a^{\centerdot}} (15,10)*+{A(CB)}; (45,10)*+{(AC)B} };
{\ar^{R1} (45,10)*+{(AC)B}; (60,0)*+{(CA)B} };
{\ar_{a^{\centerdot}} (0,0)*+{A(BC)}; (15,-10)*+{(AB)C} };
{\ar_{R} (15,-10)*+{(AB)C}; (45,-10)*+{C(AB)} };
{\ar_{a^{\centerdot}} (45,-10)*+{C(AB)}; (60,0)*+{(CA)B} };
{\ar@{=>}^{R_{(A,B|C)}} (29, 7)*{}; (29,-7)*{} }
\endxy
\]
\end{itemize}
all subject to the following four axioms.
\[
\xy
{\ar^{a} (0,0)*+{\scriptstyle (AB)(CD)}; (8,8)*+{\scriptstyle A(B(CD))} };
{\ar^{1(1R)} (8,8)*+{\scriptstyle A(B(CD))}; (16,16)*+{\scriptstyle A(B(DC))} };
{\ar^{1a^{\centerdot}} (16,16)*+{\scriptstyle A(B(DC))}; (24,24)*+{\scriptstyle A((BD)C)} };
{\ar^{1(R1)} (24,24)*+{\scriptstyle A((BD)C)}; (50,24)*+{\scriptstyle A((DB)C)} };
{\ar^{1a} (50,24)*+{\scriptstyle A((DB)C)}; (76,24)*+{\scriptstyle A(D(BC))} };
{\ar^{a^{\centerdot}} (76,24)*+{\scriptstyle A(D(BC))}; (84,16)*+{\scriptstyle (AD)(BC)} };
{\ar^{a^{\centerdot}} (84,16)*+{\scriptstyle (AD)(BC)}; (92,8)*+{\scriptstyle ((AD)B)C} };
{\ar^{(R1)1} (92,8)*+{\scriptstyle ((AD)B)C}; (100,0)*+{\scriptstyle ((DA)B)C} };
{\ar_{a^{\centerdot}} (0,0)*+{\scriptstyle (AB)(CD)}; (20,-14)*+{\scriptstyle ((AB)C)D} };
{\ar_{R} (20,-14)*+{\scriptstyle ((AB)C)D}; (50,-14)*+{\scriptstyle D((AB)C)} };
{\ar_{a^{\centerdot}} (50,-14)*+{\scriptstyle D((AB)C)}; (80,-14)*+{\scriptstyle (D(AB))C} };
{\ar_{a^{\centerdot}1} (80,-14)*+{\scriptstyle (D(AB))C}; (100,0)*+{\scriptstyle ((DA)B)C} };
{\ar^{1a^{\centerdot}} (8,8)*+{\scriptstyle A(B(CD))}; (42,16)*+{\scriptstyle A((BC)D)} };
{\ar^{1R} (42,16)*+{\scriptstyle A((BC)D)}; (76,24)*+{\scriptstyle A(D(BC))} };
{\ar_{a^{\centerdot}} (42,16)*+{\scriptstyle A((BC)D)}; (30,0)*+{\scriptstyle (A(BC))D} };
{\ar_{a^{\centerdot}1} (30,0)*+{\scriptstyle (A(BC))D}; (20,-14)*+{\scriptstyle ((AB)C)D} };
{\ar^{R} (30,0)*+{\scriptstyle (A(BC))D}; (50,0)*+{\scriptstyle D(A(BC))} };
{\ar^{a^{\centerdot}} (50,0)*+{\scriptstyle D(A(BC))}; (74,0)*+{\scriptstyle (DA)(BC)} };
{\ar_{a^{\centerdot}} (74,0)*+{\scriptstyle (DA)(BC)}; (100,0)*+{\scriptstyle ((DA)B)C} };
{\ar^{1a^{\centerdot}} (50,0)*+{\scriptstyle D(A(BC))}; (50,-14)*+{\scriptstyle D((AB)C)} };
{\ar_{R1} (84,16)*+{\scriptstyle (AD)(BC)}; (74,0)*+{\scriptstyle (DA)(BC)} };
(16,0)*{\Downarrow \pi}; (34,20)*{\Downarrow 1R_{(B,C|D)}}; (62,11)*{\Downarrow R_{(A,BC|D)}}; (36,-7)*{\cong}; (87,4)*{\cong}; (70,-6)*{\Downarrow \pi};
{\ar@{=} (50,-16)*{}; (50,-25)*{} };
{\ar^{a} (0,-51)*+{\scriptstyle (AB)(CD)}; (8,-43)*+{\scriptstyle A(B(CD))} };
{\ar^{1(1R)} (8,-43)*+{\scriptstyle A(B(CD))}; (16,-35)*+{\scriptstyle A(B(DC))} };
{\ar^{1a^{\centerdot}} (16,-35)*+{\scriptstyle A(B(DC))}; (24,-27)*+{\scriptstyle A((BD)C)} };
{\ar^{1(R1)} (24,-27)*+{\scriptstyle A((BD)C)}; (50,-27)*+{\scriptstyle A((DB)C)} };
{\ar^{1a} (50,-27)*+{\scriptstyle A((DB)C)}; (76,-27)*+{\scriptstyle A(D(BC))} };
{\ar^{a^{\centerdot}} (76,-27)*+{\scriptstyle A(D(BC))}; (84,-35)*+{\scriptstyle (AD)(BC)} };
{\ar^{a^{\centerdot}} (84,-35)*+{\scriptstyle (AD)(BC)}; (92,-43)*+{\scriptstyle ((AD)B)C} };
{\ar^{(R1)1} (92,-43)*+{\scriptstyle ((AD)B)C}; (100,-51)*+{\scriptstyle ((DA)B)C} };
{\ar_{a^{\centerdot}} (0,-51)*+{\scriptstyle (AB)(CD)}; (20,-67.5)*+{\scriptstyle ((AB)C)D} };
{\ar_{R} (20,-67.5)*+{\scriptstyle ((AB)C)D}; (50,-67.5)*+{\scriptstyle D((AB)C)} };
{\ar_{a^{\centerdot}} (50,-67.5)*+{\scriptstyle D((AB)C)}; (80,-67.5)*+{\scriptstyle (D(AB))C} };
{\ar_{a^{\centerdot}1} (80,-67.5)*+{\scriptstyle (D(AB))C}; (100,-51)*+{\scriptstyle ((DA)B)C} };
{\ar_{1R} (0,-51)*+{\scriptstyle (AB)(CD)}; (24,-51)*+{\scriptstyle (AB)(DC)} };
{\ar^{a} (24,-51)*+{\scriptstyle (AB)(DC)}; (16,-35)*+{\scriptstyle A(B(DC))} };
{\ar_{a^{\centerdot}} (24,-51)*+{\scriptstyle (AB)(DC)}; (48,-51)*+{\scriptstyle ((AB)D)C} };
{\ar_{a1} (48,-51)*+{\scriptstyle ((AB)D)C}; (48,-44)*+{\scriptstyle (A(BD))C} };
{\ar_{(1R)1} (48,-44)*+{\scriptstyle (A(BD))C}; (65,-35)*+{\scriptstyle (A(DB))C} };
{\ar^{a^{\centerdot}} (48,-44)*+{\scriptstyle (A(BD))C}; (24,-27)*+{\scriptstyle A((BD)C)} };
{\ar_{a^{\centerdot}} (65,-35)*+{\scriptstyle (A(DB))C}; (50,-27)*+{\scriptstyle A((DB)C)} };
{\ar_{a^{\centerdot}1} (65,-35)*+{\scriptstyle (A(DB))C}; (92,-43)*+{\scriptstyle ((AD)B)C} };
{\ar^{R1} (48,-51)*+{\scriptstyle ((AB)D)C}; (80,-67.5)*+{\scriptstyle (D(AB))C} };
(10,-47)*{\cong}; (31,-43)*{\Downarrow \pi}; (48,-35)*{\cong}; (70,-31)*{\Downarrow \pi}; (72,-50)*{\Downarrow R_{(A,B|D)}1}; (30,-59)*{\Downarrow R_{(AB,C|D)}}\endxy
\]

\[
\xy
{\ar^{R1} (0,0)*+{\scriptstyle (AB)(CD)}; (7,9)*+{\scriptstyle (BA)(CD)} };
{\ar^{a^{\centerdot}} (7,9)*+{\scriptstyle (BA)(CD)}; (14,18)*+{\scriptstyle ((BA)C)D} };
{\ar^{a1} (14,18)*+{\scriptstyle ((BA)C)D}; (38,18)*+{\scriptstyle (B(AC))D} };
{\ar^{(1R)1} (38,18)*+{\scriptstyle (B(AC))D}; (62,18)*+{\scriptstyle (B(CA))D} };
{\ar^{a^{\centerdot}1} (62,18)*+{\scriptstyle (B(CA))D}; (86,18)*+{\scriptstyle ((BC)A)D} };
{\ar^{a} (86,18)*+{\scriptstyle ((BC)A)D}; (93,9)*+{\scriptstyle (BC)(AD)} };
{\ar^{1R} (93,9)*+{\scriptstyle (BC)(AD)}; (100,0)*+{\scriptstyle (BC)(DA)} };
{\ar_{a} (0,0)*+{\scriptstyle (AB)(CD)}; (12,-14)*+{\scriptstyle A(B(CD))} };
{\ar_{R}  (12,-14)*+{\scriptstyle A(B(CD))}; (37,-14)*+{\scriptstyle (B(CD))A} };
{\ar_{a} (37,-14)*+{\scriptstyle (B(CD))A}; (63,-14)*+{\scriptstyle B((CD)A)} };
{\ar_{1a} (63,-14)*+{\scriptstyle B((CD)A)}; (88,-14)*+{\scriptstyle B(C(DA))} };
{\ar_{a^{\centerdot}} (88,-14)*+{\scriptstyle B(C(DA))}; (100,0)*+{\scriptstyle (BC)(DA)} };
{\ar_{a} (7,9)*+{\scriptstyle (BA)(CD)}; (25,0)*+{\scriptstyle B(A(CD))} };
{\ar_{1 a^{\centerdot}} (25,0)*+{\scriptstyle B(A(CD))}; (40,7)*+{\scriptstyle B((AC)D)} };
{\ar^{a^{\centerdot}} (40,7)*+{\scriptstyle B((AC)D)}; (38,18)*+{\scriptstyle (B(AC))D} };
{\ar_{1R} (25,0)*+{\scriptstyle B(A(CD))}; (63,-14)*+{\scriptstyle B((CD)A)} };
{\ar_{1(R1)} (40,7)*+{\scriptstyle B((AC)D)}; (62,6)*+{\scriptstyle B((CA)D)} };
{\ar^{a^{\centerdot}} (62,6)*+{\scriptstyle B((CA)D)}; (62,18)*+{\scriptstyle (B(CA))D} };
{\ar_{1a} (62,6)*+{\scriptstyle B((CA)D)}; (75,-4)*+{\scriptstyle B(C(AD))} };
{\ar_{1(1R)} (75,-4)*+{\scriptstyle B(C(AD))}; (88,-14)*+{\scriptstyle B(C(DA))} };
{\ar^{a^{\centerdot}} (75,-4)*+{\scriptstyle B(C(AD))}; (93,9)*+{\scriptstyle (BC)(AD)} };
(26,10)*{\Downarrow \pi}; (50,11)*{\cong}; (74,10)*{\Downarrow \pi}; (20,-7)*{\Downarrow R_{(A|B,CD)}}; (50,-3)*{\Downarrow 1 R_{(A|C,D)}}; (88,-5)*{\cong};
{\ar@{=} (50,-17)*{}; (50,-25)*{} };
{\ar^{R1} (0,-47)*+{\scriptstyle (AB)(CD)}; (7,-38)*+{\scriptstyle (BA)(CD)} };
{\ar^{a^{\centerdot}} (7,-38)*+{\scriptstyle (BA)(CD)}; (14,-29)*+{\scriptstyle ((BA)C)D} };
{\ar^{a1} (14,-29)*+{\scriptstyle ((BA)C)D}; (38,-29)*+{\scriptstyle (B(AC))D} };
{\ar^{(1R)1} (38,-29)*+{\scriptstyle (B(AC))D}; (62,-29)*+{\scriptstyle (B(CA))D} };
{\ar^{a^{\centerdot}1} (62,-29)*+{\scriptstyle (B(CA))D}; (86,-29)*+{\scriptstyle ((BC)A)D} };
{\ar^{a} (86,-29)*+{\scriptstyle ((BC)A)D}; (93,-38)*+{\scriptstyle (BC)(AD)} };
{\ar^{1R} (93,-38)*+{\scriptstyle (BC)(AD)}; (100,-47)*+{\scriptstyle (BC)(DA)} };
{\ar_{a} (0,-47)*+{\scriptstyle (AB)(CD)}; (12,-61)*+{\scriptstyle A(B(CD))} };
{\ar_{R}  (12,-61)*+{\scriptstyle A(B(CD))}; (37,-61)*+{\scriptstyle (B(CD))A} };
{\ar_{a} (37,-61)*+{\scriptstyle (B(CD))A}; (63,-61)*+{\scriptstyle B((CD)A)} };
{\ar_{1a} (63,-61)*+{\scriptstyle B((CD)A)}; (88,-61)*+{\scriptstyle B(C(DA))} };
{\ar_{a^{\centerdot}} (88,-61)*+{\scriptstyle B(C(DA))}; (100,-47)*+{\scriptstyle (BC)(DA)} };
{\ar_{a^{\centerdot}} (0,-47)*+{\scriptstyle (AB)(CD)}; (29,-41)*+{\scriptstyle ((AB)C)D} };
{\ar_{(R1)1} (29,-41)*+{\scriptstyle ((AB)C)D}; (14,-29)*+{\scriptstyle ((BA)C)D} };
{\ar_{a1} (29,-41)*+{\scriptstyle ((AB)C)D}; (57,-35)*+{\scriptstyle (A(BC))D} };
{\ar_{R1} (57,-35)*+{\scriptstyle (A(BC))D}; (86,-29)*+{\scriptstyle ((BC)A)D} };
{\ar^{a} (57,-35)*+{\scriptstyle (A(BC))D}; (59,-43)*+{\scriptstyle A((BC)D)} };
{\ar^{R} (59,-43)*+{\scriptstyle A((BC)D)}; (61,-52)*+{\scriptstyle ((BC)D)A} };
{\ar^{1a^{\centerdot}} (12,-61)*+{\scriptstyle A(B(CD))}; (59,-43)*+{\scriptstyle A((BC)D)} };
{\ar_{a^{\centerdot}1} (37,-61)*+{\scriptstyle (B(CD))A}; (61,-52)*+{\scriptstyle ((BC)D)A} };
{\ar_{a} (61,-52)*+{\scriptstyle ((BC)D)A}; (100,-47)*+{\scriptstyle (BC)(DA)} };
(16,-40)*{\cong}; (40,-34)*{\Downarrow R_{(A|B,C)}1}; (28,-48)*{\Downarrow \pi}; (45,-52)*{\cong}; (77,-43)*{\Downarrow R_{(A|BC,D)}}; (76,-56)*{\Downarrow \pi}
\endxy
\]

\[
\xy
{\ar^{\scriptstyle (1R)1} (0,0)*+{\scriptscriptstyle (A(BC))D}; (0,12.6)*+{\scriptscriptstyle (A(CB))D} };
{\ar^{\scriptstyle a^{\centerdot}1} (0,12.6)*+{\scriptscriptstyle (A(CB))D}; (8,25.2)*+{\scriptscriptstyle ((AC)B)D} };
{\ar^{\scriptstyle a} (8,25.2)*+{\scriptscriptstyle ((AC)B)D}; (29,25.2)*+{\scriptscriptstyle (AC)(BD)} };
{\ar^{\scriptstyle R1} (29,25.2)*+{\scriptscriptstyle (AC)(BD)}; (50,25.2)*+{\scriptscriptstyle (CA)(BD)} };
{\ar^{\scriptstyle 1R} (50,25.2)*+{\scriptscriptstyle (CA)(BD)}; (71,25.2)*+{\scriptscriptstyle (CA)(DB)} };
{\ar^{\scriptstyle a} (71,25.2)*+{\scriptscriptstyle (CA)(DB)}; (92,25.2)*+{\scriptscriptstyle C(A(DB))} };
{\ar^{\scriptstyle 1a^{\centerdot}} (92,25.2)*+{\scriptscriptstyle C(A(DB))}; (100,12.6)*+{\scriptscriptstyle C((AD)B)} };
{\ar^{\scriptstyle 1(R1)} (100,12.6)*+{\scriptscriptstyle C((AD)B)}; (100,0)*+{\scriptscriptstyle C((DA)B)} };
{\ar_{\scriptstyle a^{\centerdot}1} (0,0)*+{\scriptscriptstyle (A(BC))D}; (17,-15.75)*+{\scriptscriptstyle ((AB)C)D} };
{\ar_{\scriptstyle a} (17,-15.75)*+{\scriptscriptstyle ((AB)C)D}; (39,-15.75)*+{\scriptscriptstyle (AB)(CD)} };
{\ar_{\scriptstyle R} (39,-15.75)*+{\scriptscriptstyle (AB)(CD)}; (61,-15.75)*+{\scriptscriptstyle (CD)(AB)} };
{\ar_{\scriptstyle a} (61,-15.75)*+{\scriptscriptstyle (CD)(AB)}; (83,-15.75)*+{\scriptscriptstyle C(D(AB))} };
{\ar_{\scriptstyle 1 a^{\centerdot}} (83,-15.75)*+{\scriptscriptstyle C(D(AB))}; (100,0)*+{\scriptscriptstyle C((DA)B)} };
{\ar^{\scriptstyle R1} (17,-15.75)*+{\scriptscriptstyle ((AB)C)D}; (23,-2.1)*+{\scriptscriptstyle (C(AB))D} };
{\ar^{\scriptstyle a^{\centerdot}1} (23,-2.1)*+{\scriptscriptstyle (C(AB))D}; (29,11.55)*+{\scriptscriptstyle ((CA)B)D} };
{\ar_{\scriptstyle (R1)1} (8,25.2)*+{\scriptscriptstyle ((AC)B)D}; (29,11.55)*+{\scriptscriptstyle ((CA)B)D} };
{\ar_{\scriptstyle a} (29,11.55)*+{\scriptscriptstyle ((CA)B)D}; (50,25.2)*+{\scriptscriptstyle (CA)(BD)} };
{\ar_{\scriptstyle a} (50,25.2)*+{\scriptscriptstyle (CA)(BD)}; (71,11.55)*+{\scriptscriptstyle C(A(BD))} };
{\ar_{\scriptstyle 1(1R)} (71,11.55)*+{\scriptscriptstyle C(A(BD))}; (92,25.2)*+{\scriptscriptstyle C(A(DB))} };
{\ar^{\scriptstyle 1a^{\centerdot}} (71,11.55)*+{\scriptscriptstyle C(A(BD))}; (60,-1.05)*+{\scriptscriptstyle C((AB)D)} };
{\ar^{\scriptstyle a} (23,-2.1)*+{\scriptscriptstyle (C(AB))D}; (60,-1.05)*+{\scriptscriptstyle C((AB)D)} };
{\ar^{\scriptstyle 1R} (60,-1.05)*+{\scriptscriptstyle C((AB)D)}; (83,-15.75)*+{\scriptscriptstyle C(D(AB))} };
(12,8)*{\scriptstyle \Downarrow R_{(A,B|C)}1}; (29,18)*{\scriptstyle \cong}; (50,7)*{\scriptstyle \Downarrow \pi}; (71,18)*{\scriptstyle \cong}; (83,6)*{\scriptstyle \Downarrow 1 R_{(A,B|D)}}; (50,-8)*{\scriptstyle \Downarrow R_{(AB|C,D)}};
{\ar@{=} (50,-20.5)*{}; (50,-27.5)*{} };
{\ar^{\scriptstyle (1R)1} (0,-64)*+{\scriptscriptstyle (A(BC))D}; (0,-48)*+{\scriptscriptstyle (A(CB))D} };
{\ar^{\scriptstyle a^{\centerdot}1} (0,-48)*+{\scriptscriptstyle (A(CB))D}; (8,-30)*+{\scriptscriptstyle ((AC)B)D} };
{\ar^{\scriptstyle a} (8,-30)*+{\scriptscriptstyle ((AC)B)D}; (29,-30)*+{\scriptscriptstyle (AC)(BD)} };
{\ar^{\scriptstyle R1} (29,-30)*+{\scriptscriptstyle (AC)(BD)}; (50,-30)*+{\scriptscriptstyle (CA)(BD)} };
{\ar^{\scriptstyle 1R} (50,-30)*+{\scriptscriptstyle (CA)(BD)}; (71,-30)*+{\scriptscriptstyle (CA)(DB)} };
{\ar^{\scriptstyle a} (71,-30)*+{\scriptscriptstyle (CA)(DB)}; (92,-30)*+{\scriptscriptstyle C(A(DB))} };
{\ar^{\scriptstyle 1a^{\centerdot}} (92,-30)*+{\scriptscriptstyle C(A(DB))}; (100,-48)*+{\scriptscriptstyle C((AD)B)} };
{\ar^{\scriptstyle 1(R1)} (100,-48)*+{\scriptscriptstyle C((AD)B)}; (100,-64)*+{\scriptscriptstyle C((DA)B)} };
{\ar_{\scriptstyle a^{\centerdot}1} (0,-64)*+{\scriptscriptstyle (A(BC))D}; (17,-82)*+{\scriptscriptstyle ((AB)C)D} };
{\ar_{\scriptstyle a} (17,-82)*+{\scriptscriptstyle ((AB)C)D}; (39,-82)*+{\scriptscriptstyle (AB)(CD)} };
{\ar_{\scriptstyle R} (39,-82)*+{\scriptscriptstyle (AB)(CD)}; (61,-82)*+{\scriptscriptstyle (CD)(AB)} };
{\ar_{\scriptstyle a} (61,-82)*+{\scriptscriptstyle (CD)(AB)}; (83,-82)*+{\scriptscriptstyle C(D(AB))} };
{\ar_{\scriptstyle 1 a^{\centerdot}} (83,-82)*+{\scriptscriptstyle C(D(AB))}; (100,-64)*+{\scriptscriptstyle C((DA)B)} };
{\ar^{\scriptstyle a} (0,-64)*+{\scriptscriptstyle (A(BC))D}; (15,-57)*+{\scriptscriptstyle A((BC)D)} };
{\ar^{\scriptstyle 1(R1)} (15,-57)*+{\scriptscriptstyle A((BC)D)}; (15,-45)*+{\scriptscriptstyle A((CB)D)} };
{\ar_{\scriptstyle a} (0,-48)*+{\scriptscriptstyle (A(CB))D}; (15,-45)*+{\scriptscriptstyle A((CB)D)} };
{\ar_{\scriptstyle 1a} (15,-45)*+{\scriptscriptstyle A((CB)D)}; (26,-38)*+{\scriptscriptstyle A(C(BD))} };
{\ar_{\scriptstyle a^{\centerdot}} (26,-38)*+{\scriptscriptstyle A(C(BD))}; (29,-30)*+{\scriptscriptstyle (AC)(BD)} };
{\ar_{\scriptstyle 1R} (29,-30)*+{\scriptscriptstyle (AC)(BD)}; (50,-42)*+{\scriptscriptstyle (AC)(DB)} };
{\ar_{\scriptstyle R1} (50,-42)*+{\scriptscriptstyle (AC)(DB)}; (71,-30)*+{\scriptscriptstyle (CA)(DB)} };
{\ar_{\scriptstyle a^{\centerdot}} (71,-30)*+{\scriptscriptstyle (CA)(DB)}; (68,-49)*+{\scriptscriptstyle ((CA)D)B} };
{\ar_{\scriptstyle a1} (68,-49)*+{\scriptscriptstyle ((CA)D)B}; (84,-45)*+{\scriptscriptstyle (C(AD))B} };
{\ar_{\scriptstyle a} (84,-45)*+{\scriptscriptstyle (C(AD))B}; (100,-48)*+{\scriptscriptstyle C((AD)B)} };
{\ar^{\scriptstyle (1R)1} (84,-45)*+{\scriptscriptstyle (C(AD))B}; (84,-60)*+{\scriptscriptstyle (C(DA))B} };
{\ar_{\scriptstyle a} (84,-60)*+{\scriptscriptstyle (C(DA))B}; (100,-64)*+{\scriptscriptstyle C((DA)B)} };
{\ar_{\scriptstyle a^{\centerdot}} (61,-82)*+{\scriptscriptstyle (CD)(AB)}; (75,-71)*+{\scriptscriptstyle ((CD)A)B} };
{\ar_{\scriptstyle a1} (75,-71)*+{\scriptscriptstyle ((CD)A)B};  (84,-60)*+{\scriptscriptstyle (C(DA))B} };
{\ar_{\scriptstyle 1(1R)} (26,-38)*+{\scriptscriptstyle A(C(BD))}; (38,-50)*+{\scriptscriptstyle A(C(DB))} };
{\ar_{\scriptstyle a^{\centerdot}} (38,-50)*+{\scriptscriptstyle A(C(DB))}; (50,-42)*+{\scriptscriptstyle (AC)(DB)} };
{\ar_{\scriptstyle a^{\centerdot}} (50,-42)*+{\scriptscriptstyle (AC)(DB)}; (60,-59)*+{\scriptscriptstyle ((AC)D)B} };
{\ar_{\scriptstyle (R1)1} (60,-59)*+{\scriptscriptstyle ((AC)D)B}; (68,-49)*+{\scriptscriptstyle ((CA)D)B} };
{\ar^{\scriptstyle 1a} (15,-57)*+{\scriptscriptstyle A((BC)D)}; (24.5,-70)*+{\scriptscriptstyle A(B(CD))} };
{\ar^{\scriptstyle a^{\centerdot}} (24.5,-70)*+{\scriptscriptstyle A(B(CD))}; (39,-82)*+{\scriptscriptstyle (AB)(CD)} };
{\ar_{\scriptstyle 1R} (24.5,-70)*+{\scriptscriptstyle A(B(CD))}; (32,-59)*+{\scriptscriptstyle A((CD)B)} };
{\ar_{\scriptstyle 1a} (32,-59)*+{\scriptscriptstyle A((CD)B)}; (38,-50)*+{\scriptscriptstyle A(C(DB))} };
{\ar_{\scriptstyle a^{\centerdot}} (32,-59)*+{\scriptscriptstyle A((CD)B)}; (50,-71)*+{\scriptscriptstyle (A(CD))B} };
{\ar_{\scriptstyle R1} (50,-71)*+{\scriptscriptstyle (A(CD))B}; (75,-71)*+{\scriptscriptstyle ((CD)A)B} };
{\ar_{\scriptstyle a^{\centerdot}1} (50,-71)*+{\scriptscriptstyle (A(CD))B}; (60,-59)*+{\scriptscriptstyle ((AC)D)B} };
(7,-54)*{\scriptstyle \cong};(11,-38)*{\scriptstyle \Downarrow \pi}; (50,-35)*{\scriptstyle \cong}; (35,-42)*{\scriptstyle \cong}; (60,-45)*{\scriptstyle \cong}; (84,-38)*{\scriptstyle \Downarrow \pi}; (91,-55)*{\scriptstyle \cong}; (83,-76)*{\scriptstyle \Downarrow \pi}; (13,-68)*{\scriptstyle \Downarrow \pi}; (47,-60)*{\scriptstyle \Downarrow \pi}; (25,-53)*{\scriptstyle \Downarrow 1 R_{(B|C,D)}}; (50,-76)*{\scriptstyle \Downarrow R_{(A,B|CD)}}; (69,-64)*{\scriptstyle \Downarrow R_{(A|C,D)}1}
\endxy
\]

\[
\xy
{\ar^{R1} (0,0)*+{(AB)C}; (10,15)*+{(BA)C} };
{\ar^{a} (10,15)*+{(BA)C}; (30,15)*+{B(AC)} };
{\ar^{1R} (30,15)*+{B(AC)}; (50,15)*+{B(CA)} };
{\ar^{a^{\centerdot}} (50,15)*+{B(CA)}; (70,15)*+{(BC)A} };
{\ar^{R1} (70,15)*+{(BC)A}; (90,15)*+{(CB)A} };
{\ar^{a} (90,15)*+{(CB)A}; (100,0)*+{C(BA)} };
{\ar_{a} (0,0)*+{(AB)C}; (10,-15)*+{A(BC)} };
{\ar_{1R} (10,-15)*+{A(BC)}; (30,-15)*+{A(CB)} };
{\ar_{a^{\centerdot}} (30,-15)*+{A(CB)}; (50,-15)*+{(AC)B} };
{\ar_{R1} (50,-15)*+{(AC)B}; (70,-15)*+{(CA)B} };
{\ar_{a} (70,-15)*+{(CA)B}; (90,-15)*+{C(AB)} };
{\ar_{1R} (90,-15)*+{C(AB)}; (100,0)*+{C(BA)} };
{\ar^{R} (10,-15)*+{A(BC)}; (30,0)*+{(BC)A} };
{\ar^{a}  (30,0)*+{(BC)A}; (50,15)*+{B(CA)} };
{\ar_{R1} (30,0)*+{(BC)A}; (90,15)*+{(CB)A} };
{\ar_{R} (30,-15)*+{A(CB)}; (90,15)*+{(CB)A} };
(18,3.75)*{\scriptstyle \Downarrow R_{(A|B,C)}}; (58,10)*{\cong}; (40,-3.75)*{\cong}; (75,-3.75)*{\scriptstyle \Downarrow R_{(A|C,B)}^{-1}};
{\ar@{=} (50,-19.5)*{}; (50,-31.5)*{} };
{\ar^{R1} (0,-52.5)*+{(AB)C}; (10,-37.5)*+{(BA)C} };
{\ar^{a} (10,-37.5)*+{(BA)C}; (30,-37.5)*+{B(AC)} };
{\ar^{1R} (30,-37.5)*+{B(AC)}; (50,-37.5)*+{B(CA)} };
{\ar^{a^{\centerdot}} (50,-37.5)*+{B(CA)}; (70,-37.5)*+{(BC)A} };
{\ar^{R1} (70,-37.5)*+{(BC)A}; (90,-37.5)*+{(CB)A} };
{\ar^{a} (90,-37.5)*+{(CB)A}; (100,-52.5)*+{C(BA)} };
{\ar_{a} (0,-52.5)*+{(AB)C}; (10,-67.5)*+{A(BC)} };
{\ar_{1R} (10,-67.5)*+{A(BC)}; (30,-67.5)*+{A(CB)} };
{\ar_{a^{\centerdot}} (30,-67.5)*+{A(CB)}; (50,-67.5)*+{(AC)B} };
{\ar_{R1} (50,-67.5)*+{(AC)B}; (70,-67.5)*+{(CA)B} };
{\ar_{a} (70,-67.5)*+{(CA)B}; (90,-67.5)*+{C(AB)} };
{\ar_{1R} (90,-67.5)*+{C(AB)}; (100,-52.5)*+{C(BA)} };
{\ar_{R} (0,-52.5)*+{(AB)C}; (90,-67.5)*+{C(AB)} };
{\ar^{a^{\centerdot}} (30,-37.5)*+{B(AC)}; (65,-52.5)*+{(BA)C} };
{\ar_{R} (65,-52.5)*+{(BA)C}; (100,-52.5)*+{C(BA)} };
{\ar^{R1} (0,-52.5)*+{(AB)C}; (65,-52.5)*+{(BA)C} };
(20,-45)*{\cong}; (17,-61.25)*{\scriptstyle \Downarrow R_{(A,B|C)}^{-1}}; (70,-45)*{\scriptstyle \Downarrow R_{(B,A|C)}}; (70,-57.5)*{\cong}
\endxy
\]
A \textit{braided monoidal bicategory} is a monoidal bicategory equipped with a braiding.
\\

While these axioms might look quite daunting, they are in fact just algebraic expressions of the notion that a pair of homotopies that each start at the braid $\gamma$ and end at the braid $\gamma'$ are in fact themselves homotopic.  For instance, the first axiom concerns the case of a braid with four strands in which the first three strands are braided past the final one.  This can be done either all at once (this is the 1-cell target of the pasting diagram) or it can be done step-by-step in which strand 3 is braided past strand 4, then strand 2 is braided past strand 4, and finally strand 1 is braided past strand 4 (this is the 1-cell source).  The two different composite 2-cells which are claimed to be equal in this axiom are just two different ways to transform the step-by-step method into the all-at-once method using the algebra available in a braided monoidal bicategory.  The other three axioms also have similar interpretations.  In fact, one could take the presentation of these axioms seriously, and view them as certain three-dimensional polytopes in which the two-dimensional faces are precisely the 2-cells in each equation above.  Doing so produces polytopes discovered by Bar-Natan in \cite{bn}.

Recall that the category of 2-categories and strict 2-functors has a variety of tensor products, two of which shall be important here.  The first is the Cartesian structure, and this tensor product has the property that $- \times B$ is left adjoint to $[B, -]$ where $[B,C]$ is the 2-category of 2-functors, 2-natural transformations, and modifications.  The second is the Gray tensor product $\otimes$ which has the property that $- \otimes B$ is left adjoint to $\textbf{Hom}(B,-)$ where $\textbf{Hom}(B,C)$ is the 2-category of 2-functors, pseudonatural transformations, and modifications.  The Gray tensor product gives the category of 2-categories the structure of a closed symmetric monoidal category, and monoids in this symmetric monoidal category are called monoidal 2-categories or $\textbf{Gray}$-monoids.  One form of the coherence theorem for monoidal bicategories states that every monoidal bicategory is monoidally biequivalent to a monoidal 2-category.  In the final section of this paper, we will extend this result to the braided case and show that every braided monoidal bicategory is braided monoidally biequivalent to a braided monoidal 2-category as defined below.
\\ \\
 \textbf{Definition}  A braided monoidal bicategory $B$ is a \textit{braided monoidal 2-category} if the following conditions hold.
\begin{enumerate}
\item The underlying monoidal bicategory of $B$ is a $\mb{Gray}$-monoid.
\item The following unit conditions hold.
\begin{itemize}
\item The adjoint equivalence $\mb{R}_{I,A}$ between $I \otimes A$ and $A \otimes I$ is the identity adjoint equivalence on $A$.
\item The adjoint equivalence $\mb{R}_{A,I}$ between $A \otimes I$ and $I \otimes A$ is the identity adjoint equivalence on $A$.
\item The isomorphism 2-cell
\[
\xy
{\ar^{R} (0,0)*+{ABI}; (40,0)*+{BIA} };
{\ar_{R1} (0,0)*+{ABI}; (20,-12)*+{BAI} };
{\ar_{1R} (20,-12)*+{BAI}; (40,0)*+{BIA} };
{\ar@{=>}^{R_{(A|B,I)}} (17,-2)*{}; (17,-8)*{} }
\endxy
\]
is the identity 2-cell $1_{R_{A,B}}$.
\item Similarly, the isomorphism 2-cells $R_{(A|I,B)}, R_{(A,I|B)}, R_{(I,A|B)}$ all equal the identity 2-cell $1_{R_{A,B}}$.
\item The isomorphism 2-cells $R_{(I|A,B)}, R_{(A,B|I)}$ both equal the identity 2-cell $1_{1_{AB}}$.
\end{itemize}
\end{enumerate}

It is worth noting that in a braided monoidal 2-category, not only are there additional unit axioms, but the standard four axioms for a braiding are substantially simpler.  For instance, the fourth axiom becomes the equality of pastings below.
\[
\xy
(0,0)*+{ABC}="1"; (10,12)*+{BAC}="2"; (35,12)*+{BCA}="3"; (45,0)*+{CBA}="4";
(10,-12)*+{ACB}="5"; (35,-12)*+{CAB}="6";
(62.5,0)*+{ABC}="7"; (72.5,12)*+{BAC}="8"; (97.5,12)*+{BCA}="9"; (107.5,0)*+{CBA}="10";
(72.5,-12)*+{ACB}="11"; (97.5,-12)*+{CAB}="12";
{\ar^{R1} "1"; "2"};
{\ar^{1R} "2"; "3"};
{\ar^{R1} "3"; "4"};
{\ar_{1R} "1"; "5"};
{\ar_{R1} "5"; "6"};
{\ar_{1R} "6"; "4"};
{\ar_{R} "1"; "3"};
{\ar^{R} "5"; "4"};
{\ar^{R1} "7"; "8"};
{\ar^{1R} "8"; "9"};
{\ar^{R1} "9"; "10"};
{\ar_{1R} "7"; "11"};
{\ar_{R1} "11"; "12"};
{\ar_{1R} "12"; "10"};
{\ar_{R} "8"; "10"};
{\ar^{R} "7"; "12"};
(53.75,0)*{=};
(22.5,0)*{\cong}; (85,0)*{\cong}; (12,8)*{\scriptstyle \Downarrow R}; (33,-8)*{\scriptstyle \Downarrow R^{-1}}; (75.5,-8)*{\scriptstyle \Downarrow R^{-1}}; (95.5,8)*{\scriptstyle \Downarrow R}
\endxy
\]

We will also need the notion of a braided monoidal functor between braided monoidal bicategories.
\\ \\
\textbf{Definition}  Let $B,C$ be braided monoidal bicategories.  Let $B,C$ be braided monoidal bicategories.  A \textit{braided monoidal functor} $F:B \rightarrow C$ consists of
\begin{itemize}
\item an underlying monoidal functor $F: B \rightarrow C$ and
\item an invertible modification $U$ as displayed below,
\[
\xy
{\ar^{F \times F} (0,0)*+{B \times B}; (30,0)*+{C \times C} };
{\ar^{\otimes_{C}} (30,0)*+{C \times C}; (30,-20)*+{C} };
{\ar_{\tau} (0,0)*+{B \times B}; (-10,-10)*+{B \times B} };
{\ar_{\otimes} (-10,-10)*+{B \times B}; (0,-20)*+{B} };
{\ar_{F} (0,-20)*+{B}; (30,-20)*+{C} };
{\ar^{\otimes} (0,0)*+{B \times B}; (0,-20)*+{B} };
(15,-10)*{\Downarrow \chi}; (-2.5,-10)*{\stackrel{R}{\Leftarrow}};
{\ar^{F \times F} (60,0)*+{B \times B}; (90,0)*+{C \times C} };
{\ar^{\otimes_{C}} (90,0)*+{C \times C}; (90,-20)*+{C} };
{\ar_{\tau} (60,0)*+{B \times B}; (50,-10)*+{B \times B} };
{\ar_{\otimes} (50,-10)*+{B \times B}; (60,-20)*+{B} };
{\ar_{F} (60,-20)*+{B}; (90,-20)*+{C} };
{\ar_{\tau} (90,0)*+{C \times C}; (80,-10)*+{C \times C} };
{\ar_{\otimes} (80,-10)*+{C \times C}; (90,-20)*+{C} };
{\ar^{F \times F} (50,-10)*+{B \times B}; (80,-10)*+{C \times C} };
(68,-5)*{=}; (68,-15)*{\Downarrow \chi}; (87.5,-10)*{\stackrel{R}{\Leftarrow}};
{\ar@{=>}^{U} (37,-10)*{}; (43,-10)*{} }
\endxy
\]
\end{itemize}
subject to the following two axioms.
\[
\xy
{\ar^{a^{\centerdot}} (0,0)*+{Fx(FyFz)}; (10,15)*+{(FxFy)Fz} };
{\ar^{R} (10,15)*+{(FxFy)Fz}; (35,25)*+{Fz(FxFy)} };
{\ar^{a^{\centerdot}} (35,25)*+{Fz(FxFy)}; (65,25)*+{(FzFx)Fy} };
{\ar^{\chi 1} (65,25)*+{(FzFx)Fy}; (90,15)*+{F(zx)Fy} };
{\ar^{ \chi} (90,15)*+{F(zx)Fy}; (100,0)*+{F((zx)y)} };
{\ar_{1 \chi} (0,0)*+{Fx(FyFz)}; (10,-15)*+{FxF(yz)} };
{\ar_{1 FR} (10,-15)*+{FxF(yz)}; (35,-25)*+{FxF(zy)} };
{\ar_{\chi} (35,-25)*+{FxF(zy)}; (65,-25)*+{F(x(zy))} };
{\ar_{Fa^{\centerdot}} (65,-25)*+{F(x(zy))}; (90,-15)*+{F((xz)y)} };
{\ar_{F(R1)} (90,-15)*+{F((xz)y)}; (100,0)*+{F((zx)y)} };
{\ar_{1R} (0,0)*+{Fx(FyFz)}; (25,0)*+{Fx(FzFy)} };
{\ar_{a^{\centerdot}} (25,0)*+{Fx(FzFy)}; (40,12.5)*+{(FxFz)Fy} };
{\ar^{R1} (40,12.5)*+{(FxFz)Fy}; (65,25)*+{(FzFx)Fy} };
{\ar^{1 \chi} (25,0)*+{Fx(FzFy)}; (35,-25)*+{FxF(zy)} };
{\ar_{\chi 1} (40,12.5)*+{(FxFz)Fy}; (65, -1.25)*+{F(xz)Fy} };
{\ar_{\chi} (65, -1.25)*+{F(xz)Fy}; (90,-15)*+{F((xz)y)} };
{\ar_{FR 1} (65, -1.25)*+{F(xz)Fy}; (90,15)*+{F(zx)Fy} };
(26,14)*+{\scriptstyle \Downarrow R}; (14,-8)*{\scriptstyle \Downarrow 1 U}; (50,-10)*{\scriptstyle \Downarrow \omega}; (65,15)*{\scriptstyle \Downarrow U 1}; (88,0)*{\cong};
{\ar^{a^{\centerdot}} (0,-60)*+{Fx(FyFz)}; (10,-45)*+{(FxFy)Fz} };
{\ar^{R} (10,-45)*+{(FxFy)Fz}; (35,-35)*+{Fz(FxFy)} };
{\ar^{a^{\centerdot}} (35,-35)*+{Fz(FxFy)}; (65,-35)*+{(FzFx)Fy} };
{\ar^{\chi 1} (65,-35)*+{(FzFx)Fy}; (90,-45)*+{F(zx)Fy} };
{\ar^{ \chi} (90,-45)*+{F(zx)Fy}; (100,-60)*+{F((zx)y)} };
{\ar_{1 \chi} (0,-60)*+{Fx(FyFz)}; (10,-75)*+{FxF(yz)} };
{\ar_{1 FR} (10,-75)*+{FxF(yz)}; (35,-85)*+{FxF(zy)} };
{\ar_{\chi} (35,-85)*+{FxF(zy)}; (65,-85)*+{F(x(zy))} };
{\ar_{Fa^{\centerdot}} (65,-85)*+{F(x(zy))}; (90,-75)*+{F((xz)y)} };
{\ar_{F(R1)} (90,-75)*+{F((xz)y)}; (100,-60)*+{F((zx)y)} };
{\ar^{\chi} (10,-75)*+{FxF(yz)}; (40,-75)*+{F(x(yz))} };
{\ar^{F(1R)} (40,-75)*+{F(x(yz))}; (65,-85)*+{F(x(zy))} };
{\ar^{Fa^{\centerdot}} (40,-75)*+{F(x(yz))}; (40,-60)*+{F((xy)z)} };
{\ar_{FR} (40,-60)*+{F((xy)z)}; (70,-60)*+{F(z(xy))} };
{\ar_{Fa^{\centerdot}} (70,-60)*+{F(z(xy))}; (100,-60)*+{F((zx)y)} };
{\ar_{\chi 1} (10,-45)*+{(FxFy)Fz}; (25,-52.5)*+{F(xy)Fz} };
{\ar_{\chi} (25,-52.5)*+{F(xy)Fz}; (40,-60)*+{F((xy)z)} };
{\ar^{1 \chi} (35,-35)*+{Fz(FxFy)}; (52.5,-47.5)*+{FzF(xy)} };
{\ar^{\chi} (52.5,-47.5)*+{FzF(xy)}; (70,-60)*+{F(z(xy))} };
{\ar_{R} (25,-52.5)*+{F(xy)Fz}; (52.5,-47.5)*+{FzF(xy)} };
(24,-60)*{\scriptstyle \Downarrow \omega}; (32,-43)*{\cong}; (50,-55)*{\scriptstyle \Downarrow U}; (72,-47)*{\scriptstyle \Downarrow \omega}; (38,-80)*{\cong}; (70,-72)*{\scriptstyle \Downarrow FR};
{\ar@{=} (50,-28)*{}; (50,-32)*{} }
\endxy
\]
\[
\xy
{\ar^{a} (0,0)*+{(FxFy)Fz}; (10,15)*+{Fx(FyFz)} };
{\ar^{R} (10,15)*+{Fx(FyFz)}; (35,25)*+{(FyFz)Fx} };
{\ar^{a} (35,25)*+{(FyFz)Fx}; (65,25)*+{Fy(FzFx)} };
{\ar^{1 \chi} (65,25)*+{Fy(FzFx)}; (90,15)*+{FyF(zx)} };
{\ar^{ \chi} (90,15)*+{FyF(zx)}; (100,0)*+{F(y(zx))} };
{\ar_{\chi 1} (0,0)*+{(FxFy)Fz}; (10,-15)*+{F(xy)Fz} };
{\ar_{ FR 1} (10,-15)*+{F(xy)Fz}; (35,-25)*+{F(yx)Fz} };
{\ar_{\chi} (35,-25)*+{F(yx)Fz}; (65,-25)*+{F((yx)z)} };
{\ar_{Fa} (65,-25)*+{F((yx)z)}; (90,-15)*+{F(y(xz))} };
{\ar_{F(1R)} (90,-15)*+{F(y(xz))}; (100,0)*+{F(y(zx))} };
{\ar_{R1} (0,0)*+{(FxFy)Fz}; (25,0)*+{(FyFx)Fz} };
{\ar_{a} (25,0)*+{(FyFx)Fz}; (40,12.5)*+{Fy(FxFz)} };
{\ar^{1R} (40,12.5)*+{Fy(FxFz)}; (65,25)*+{Fy(FzFx)} };
{\ar^{\chi 1} (25,0)*+{(FyFx)Fz}; (35,-25)*+{F(yx)Fz} };
{\ar_{1 \chi} (40,12.5)*+{Fy(FxFz)}; (65, -1.25)*+{FyF(xz)} };
{\ar_{\chi} (65, -1.25)*+{FyF(xz)}; (90,-15)*+{F(y(xz))} };
{\ar_{1FR} (65, -1.25)*+{FyF(xz)}; (90,15)*+{FyF(zx)} };
(26,14)*+{\scriptstyle \Downarrow R}; (14,-8)*{\scriptstyle \Downarrow U1}; (50,-10)*{\scriptstyle \Downarrow \omega}; (65,15)*{\scriptstyle \Downarrow 1U}; (88,0)*{\cong};
{\ar^{a} (0,-60)*+{(FxFy)Fz}; (10,-45)*+{Fx(FyFz)} };
{\ar^{R} (10,-45)*+{Fx(FyFz)}; (35,-35)*+{(FyFz)Fx} };
{\ar^{a} (35,-35)*+{(FyFz)Fx}; (65,-35)*+{Fy(FzFx)} };
{\ar^{1 \chi} (65,-35)*+{Fy(FzFx)}; (90,-45)*+{FyF(zx)} };
{\ar^{ \chi} (90,-45)*+{FyF(zx)}; (100,-60)*+{F(y(zx))} };
{\ar_{\chi 1} (0,-60)*+{(FxFy)Fz}; (10,-75)*+{F(xy)Fz} };
{\ar_{ FR 1} (10,-75)*+{F(xy)Fz}; (35,-85)*+{F(yx)Fz} };
{\ar_{\chi} (35,-85)*+{F(yx)Fz}; (65,-85)*+{F((yx)z)} };
{\ar_{Fa} (65,-85)*+{F((yx)z)}; (90,-75)*+{F(y(xz))} };
{\ar_{F(1R)} (90,-75)*+{F(y(xz))}; (100,-60)*+{F(y(zx))} };
{\ar^{\chi} (10,-75)*+{F(xy)Fz}; (40,-75)*+{F((xy)z)} };
{\ar^{F(R1)} (40,-75)*+{F((xy)z)}; (65,-85)*+{F((yx)z)} };
{\ar^{Fa} (40,-75)*+{F((xy)z)}; (40,-60)*+{F(x(yz))} };
{\ar_{FR} (40,-60)*+{F(x(yz))}; (70,-60)*+{F((yz)x)} };
{\ar_{Fa} (70,-60)*+{F((yz)x)}; (100,-60)*+{F(y(zx))} };
{\ar_{1 \chi} (10,-45)*+{Fx(FyFz)}; (25,-52.5)*+{FxF(yz)} };
{\ar_{\chi} (25,-52.5)*+{FxF(yz)}; (40,-60)*+{F(x(yz))} };
{\ar^{\chi 1} (35,-35)*+{(FyFz)Fx}; (52.5,-47.5)*+{F(yz)Fx} };
{\ar^{\chi} (52.5,-47.5)*+{F(yz)Fx}; (70,-60)*+{F((yz)x)} };
{\ar_{R} (25,-52.5)*+{FxF(yz)}; (52.5,-47.5)*+{F(yz)Fx} };
(24,-60)*{\scriptstyle \Downarrow \omega}; (32,-43)*{\cong}; (50,-55)*{\scriptstyle \Downarrow U}; (72,-47)*{\scriptstyle \Downarrow \omega}; (38,-80)*{\cong}; (70,-72)*{\scriptstyle \Downarrow FR};
{\ar@{=} (50,-28)*{}; (50,-32)*{} }
\endxy
\]

\noindent \textbf{Definition}  A functor $F:B \rightarrow C$ between braided monoidal bicategories is a \textit{braided monoidal biequivalence} if it is a braided monoidal functor and a biequivalence on the underlying bicategories.

\section{Fundamental 2-groupoids}
This section focuses on the construction of the fundamental 2-groupoid of a space.  The final goal of this section is to equip the fundamental 2-groupoid of a $\mathcal{C}_{2}$-algebra with the structure of a braided monoidal bicategory which we will then use to study the fundamental 2-groupoid of configuration spaces.

To achieve this goal, we proceed in several steps.  We first construct a tricategory of topological spaces, $\mb{Top}$.  That such a tricategory should exist is well-known, but the author knows of no reference giving an explicit construction.  The second step is to describe the fundamental 2-groupoid $\Pi_{2}$ as a functor of tricategories
\[
\Pi_{2}: \mb{Top} \rightarrow \mb{Bicat}.
\]
Our functor $\Pi_{2}$ agrees with the construction in \cite{hkk} (there called the homotopy bigroupoid), but they do not investigate the action on higher cells.  Finally we show that $\Pi_{2}X$ is a monoidal (resp., braided monoidal) bicategory when $X$ is an algebra for the operad $\mathcal{C}_{1}$ (resp., $\mathcal{C}_{2}$).

$\mb{Notation.}$  For the rest of the paper, any graphical representations of homotopies are to be read from the bottom to the top.  In the case of maps $I \times I \rightarrow X$, the first copy of $I$ will have coordinate $s$ along the horizontal axis and the second copy will have coordinate $t$ along the vertical axis.

\subsection{The tricategory $\mb{Top}$}
We begin by constructing a tricategory of topological spaces with the following cells:
\begin{itemize}
\item 0-cells are spaces,
\item 1-cells are continuous maps,
\item 2-cells are homotopies between continuous maps, and
\item 3-cells are homotopy classes of homotopies between 2-cells.
\end{itemize}
It should be noted that there is an obvious pointed version of this tricategory whose cells are based spaces, based maps, based homotopies, and based homotopy classes of homotopies between those.

\begin{prop}
Let $X,Y$ be spaces.  Then there is a bicategory  $\mb{Top}(X,Y)$ with
\begin{itemize}
\item 0-cells the continuous maps $f:X \rightarrow Y$,
\item 1-cells the homotopies $H:X \times I \rightarrow Y$ from $f$ to $g$, and
\item 2-cells the equivalence classes $[\alpha]$ of homotopies $\alpha:X \times I \times I \rightarrow Y$ such that at each time $t$, $\alpha(-,-,t)$ is a homotopy $f \Rightarrow g$, and $\alpha \sim \beta$ if there exists a homotopy $\Gamma:X \times I \times I \times I \rightarrow Y$ such that at each time $u$, $\Gamma(-,-,-,u)$ is a homotopy $H \Rrightarrow J$ with the property that $\Gamma(-,-,t,u)$ is a homotopy $f \Rightarrow g$.
\end{itemize}
\end{prop}
$\mb{Proof.}$
First, it is clear that the 1- and 2-cells with fixed source and target 0-cells form a category -- the 2-cell composition of equivalence classes of homotopies is given by the standard formula
\[
\beta \circ \alpha (-,-,t) = \left\{ \begin{array}{rl}
\alpha(-,-,2t) & \quad 0 \leq t \leq 1/2 \\
\beta(-,-,2t - 1) & \quad 1/2 \leq t \leq 1 \end{array} \right.
\]
and is well-defined on equivalence classes, and is associative and unital.  The composite $H_{2} \circ H_{1}$ of 1-cells is also defined by the usual formula for composing homotopies.
 \[
H_{2} \circ H_{1} (-,s) = \left\{ \begin{array}{rl}
H_{1}(-,2s) & \quad 0 \leq s \leq 1/2 \\
H_{2}(-,2s - 1) & \quad 1/2 \leq s \leq 1 \end{array} \right.
\]
It is easy to check that this satisfies the condition to be a 1-cell.  Similarly the identity homotopy $\textrm{id}$ is easily seen to satisfy the condition to be a 1-cell.  We define the horizontal composition $\beta * \alpha$ by the formula below.
\[
\beta * \alpha(-,s,t) = \left\{ \begin{array}{rl}
\alpha(-,2s,t) & \quad 0 \leq s \leq 1/2 \\
\beta(-,2s - 1,t) & \quad 1/2 \leq s \leq 1 \end{array} \right.
\]
This is easily checked to be well-defined on equivalence classes, and to give horizontal composition the structure of a functor.  The standard reparametrization formulas giving homotopies
\[
\begin{array}{c}
a: (H \circ J) \circ K \simeq H \circ (J \circ K) \\
l: \textrm{id} \circ H \simeq H \\
r: H \circ \textrm{id} \simeq H
\end{array}
\]
satisfy the conditions to be 2-cells.  The bicategory axioms are then trivial.  $\Box$

\begin{prop}
Let $X,Y,Z$ be spaces.  Then there is a strict functor
\[
\otimes: \mb{Top}(Y,Z) \times \mb{Top}(X,Y) \rightarrow \mb{Top}(X,Z)
\]
whose value at $(g,f)$ is the composite map $gf$.
\end{prop}
$\mb{Proof.}$
We need only provide the value at pairs of 1- and 2-cells and then check that it defines a strict functor.  Given 1-cells $H:f \Rightarrow f'$ and $J:g \rightarrow g'$, we define $J \otimes H:gf \Rightarrow g'f'$ as the composite
\[
X \times I \stackrel{1 \times \Delta}{\longrightarrow} X \times I \times I \stackrel{H \times 1}{\longrightarrow} Y \times I \stackrel{J}{\longrightarrow} Z,
\]
where $\Delta$ is the diagonal map.  Written as a formula, we have
\[
(J \otimes H)(x,s) = J \Big( H(x,s), s \Big).
\]

Now given 2-cells $[\alpha]:H \Rrightarrow H'$ and $[\beta]:J \Rrightarrow J'$, we define $[\beta] \otimes [\alpha]$ as the class of the map $\beta \otimes \alpha$ defined by
\[
\beta \otimes \alpha (x,s,t) = \beta \Big( \alpha(x,s,t), s, t \Big).
\]
It is now clear that this is well-defined on equivalence classes and so constitutes a 2-cell of the target.  Furthermore, it is simple to check that this preserves composition and units at the 2-cell level, so is a functor on hom-categories.

Now we must give unit and composition constraints for this functor.  The unit 1-cell is the identity homotopy, and it is simple to check that
\[
\textrm{id}_{g} \otimes \textrm{id}_{f} = \textrm{id}_{gf},
\]
so we set the unit constraint equal to the identity homotopy as well.  For the composition constraint, note that both of the 2-cells $(J_{2} \circ J_{1}) \otimes (H_{2} \circ H_{1})$ and $(J_{2} \otimes H_{2}) \circ (J_{1} \otimes H_{1})$ are given by the formula
\[
\begin{array}{cl}
J_{1} \Big( H_{1}(x, 2s), 2s \Big) & \quad 0 \leq s \leq \frac{1}{2} \\
J_{2} \Big( H_{2}(x, 2s-1),  2s-1 \Big) & \quad \frac{1}{2} \leq s \leq 1.
\end{array}
\]
Thus we define the composition constraint to be the identity as well.   It is an easy check to verify that these assignments strictly preserve the isomorphisms $a, l, r$, thus we have given a strict functor.  $\Box$

\begin{thm}
There is a tricategory $\mb{Top}$ with objects spaces, hom-bicategories given by $\mb{Top}(X,Y)$, and composition given by the functor $\otimes$.
\end{thm}
$\mb{Proof.}$
It remains to provide the unit functor, the associativity and unit adjoint equivalences, four invertible modifications, and to check the three tricategory axioms.

The unit functor takes the value of $1_{X}$ on the unique object, the identity homotopy on the unique 1-cell, and the class of the identity homotopy on the unique 2-cell.  We then define the unit functor $* \rightarrow \mb{Top}(X,X)$ to be the strict functor taking these values since the composite of two identity homotopies is again the identity homotopy.

The associativity and left and right unit adjoint equivalences are all defined to be the identity adjoint equivalences since $\otimes$ is strictly associative and unital on 1-cells.  Since the composite of identity homotopies is still the identity, we can define all four invertible modifications of the tricategory to be identities as well.  The tricategory axioms then follow trivially.  $\Box$

$\mb{Remark.}$  It should be clear that all of the results in this section should have corresponding $n$-dimensional analogues for $n>3$ including $n = \omega$.  Thus it should be possible to construct an $\omega$-category $\mb{Top}_{\omega}$; the tricategory $\mb{Top}$ should then be an appropriate 3-dimensional quotient.  Such constructions would then be relevant to studying the geometric nature of coherence for various kinds of monoidal $n$-dimensional categories for $n>2$.

\subsection{The functor $\Pi_2$}
This section develops the fundamental 2-groupoid construction as a functor of tricategories
\[
\Pi_{2}:\mb{Top} \rightarrow \mb{Bicat}.
\]
We begin by defining the action of $\Pi_{2}$ on cells.  Let $X$ be a space.  Then $\Pi_{2}X$ is the following bigroupoid (see \cite{hkk} for additional details).  The objects of $\Pi_{2}X$ are the points of $X$.  The 1-cells of $\Pi_{2}X$ from $x$ to $y$ are the paths $f:I \rightarrow X$ with $f(0) = x$ and $f(1) = y$.  The 2-cells of $\Pi_{2}X$ from $f$ to $g$ are homotopy classes of homotopies $\alpha:I \times I \rightarrow X$ with $\alpha(0,-) = f$ and $\alpha(1,-) = g$ which fix the boundary so that $\alpha(-,0) = x$ and $\alpha(-,1) = y$.  The composition $g \circ f$ of 1-cells is given by the composite
\[
I \stackrel{\times 2}{\longrightarrow} [0,2] \stackrel{f+g}{\longrightarrow} X
\]
of the multiplication by 2 map and then the sum $f+g$ which is $f(s)$ when restricted to $[0,1]$ and $g(s-1)$ when restricted to $[1,2]$.  The unit 1-cell is the constant path.  Both vertical and horizontal compositions of 2-cells are defined in the obvious fashion.  The associativity isomorphism $a:(hg)f \Rightarrow h(gf)$ is given by the class of the homotopy pictured below.

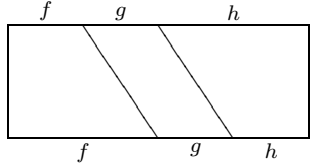
\[
\xy
{\ar@{-}_{f} (0,0)*{}; (20,0)*{} };
{\ar@{-}_{g} (20,0)*{}; (30,0)*{} };
{\ar@{-}_{h} (30,0)*{}; (40,0)*{} };
{\ar@{-}^{f} (0,15)*{}; (10,15)*{} };
{\ar@{-}^{g} (10,15)*{}; (20,15)*{} };
{\ar@{-}^{h} (20,15)*{}; (40,15)*{} };
{\ar@{-} (0,0)*{}; (0,15)*{} };
{\ar@{-} (20,0)*{}; (10,15)*{} };
{\ar@{-} (30,0)*{}; (20,15)*{} };
{\ar@{-} (40,0)*{}; (40,15)*{} };
\endxy
\]
Similarly, the left unit isomorphism $l: 1 \circ f \Rightarrow f$ is the class of the homotopy pictured below.
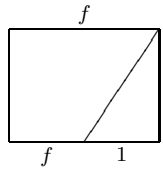
\[
\xy
{\ar@{-}_{f} (0,0)*{}; (10,0)*{} };
{\ar@{-}_{1} (10,0)*{}; (20,0)*{} };
{\ar@{-}^{f} (0,15)*{}; (20,15)*{} };
{\ar@{-} (0,0)*{}; (0,15)*{} };
{\ar@{-} (10,0)*{}; (20,15)*{} };
{\ar@{-} (20,0)*{}; (20,15)*{} };
\endxy
\]
The right unit isomorphism is defined analogously.  It is then simple to check the two bicategory axioms, and that this bicategory is in fact a bigroupoid.

Now we define $\Pi_{2}f$ for a continuous map $f:X \rightarrow Y$.  First, every continuous map sends points to points, paths to paths, and homotopies between paths to homotopies between paths; thus the action of $\Pi_{2}f$ on the cells of $\Pi_{2}X$ is obvious.  Furthermore, every continuous map takes constant paths to constant paths, thus  $\Pi_{2}f$ strictly preserves identities so we define the unit constraint to be the identity.  Similarly it is clear that $\Pi_{2}f$ strictly preserves composites, so we can define the constraint for composition to be the identity.  It is easy to check that $\Pi_{2}f$ sends the associativity and unit constraints for $\Pi_{2}X$ to those of $\Pi_{2}Y$, so $\Pi_{2}f$ becomes a strict functor.

Now we define $\Pi_{2}H$ for a homotopy $H: f \Rightarrow g$.  The component at the object $x$ is the 1-cell given by restricting $H:X \times I \rightarrow Y$ to $\{ x \} \times I$; this is a path in $Y$ which begins at $H(x,0) = f(x)$ and ends at $H(x,1) = g(x)$.  Thus we have produced a 1-cell
\[
\Pi_{2}H_{x}:\Pi_{2}f(x) \rightarrow \Pi_{2}g(x)
\]
in $\Pi_{2}Y$.  Now given a 1-cell $\gamma:x \rightarrow x'$ in $X$, we must produce a 2-cell
\[
\Pi_{2}H_{\gamma}: \Pi_{2}H_{y} \circ \Pi_{2}f(\gamma) \Rightarrow \Pi_{2}g(\gamma) \circ \Pi_{2}H_{x}
\]
in $\Pi_{2}Y$.  If we write the path giving $\Pi_{2}H_{x}$ as $H_{x}$, we are required to produce a homotopy class of homotopies
\[
H_{y} \circ f(\gamma) \Rightarrow g(\gamma) \circ H_{x}.
\]
Now the composite $H \circ (\gamma \times 1)$ is a map $I \times I \rightarrow X$ as pictured below.
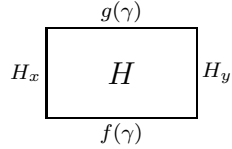
\[
\xy
{\ar@{-}^{g(\gamma)} (0,0)*+{}; (20,0)*+{} };
{\ar@{-}^{H_{y}} (20,0)*+{}; (20,-12)*+{} };
{\ar@{-}_{H_{x}} (0,0)*+{}; (0,-12)*+{} };
{\ar@{-}_{f(\gamma)} (0,-12)*+{}; (20,-12)*+{} };
(10,-6)*{H}
\endxy
\]
This is visibly not a 2-cell in $\Pi_{2}Y$ since the vertical boundaries are not constant, thus we modify this map as follows to define $\Pi_{2}H_{\gamma}$.
\[
\Pi_{2}H_{\gamma} = \left\{ \begin{array}{cl}
f \Big( \gamma(2s) \Big) & \quad s+ t \leq \frac{1}{2} \\
g \Big( \gamma(2s-1) \Big) & \quad s+ t \geq \frac{3}{2} \\
H_{x}(2s) & \quad t-s \geq \frac{1}{2} \\
H_{y}(2s-1) & \quad s-t \leq \frac{1}{2} \\
H \circ (\gamma \times 1) \Big(s-t+\frac{1}{2}, s+t - \frac{1}{2} \Big) & \quad \textrm{otherwise} \end{array} \right.
\]
We picture this map as the square below, where the corner regions with dotted lines are vertically constant and the square in the middle is $H \circ (\gamma \times 1)$ that has been shrunk and rotated.
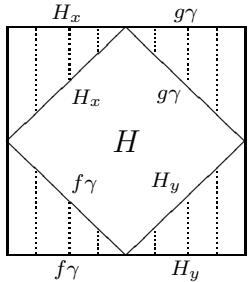
\[
\xy0;/r.15pc/:
{\ar@{-}_{f \gamma} (0,0)*{}; (25,0)*{} };
{\ar@{-}_{H_{y}} (25,0)*{}; (50,0)*{} };
{\ar@{-} (0,0)*{}; (0,48)*{} };
{\ar@{-}^{H_{x}} (0,48)*{}; (25,48)*{} };
{\ar@{-}^{g \gamma} (25,48)*{}; (50,48)*{} };
{\ar@{-} (50,0)*{}; (50,48)*{} };
{\ar@{-}_{f \gamma} (25,0)*{}; (0,24)*{} };
{\ar@{-}^{H_{y}} (25,0)*{}; (50,24)*{} };
{\ar@{-}^{H_{x}} (25,48)*{}; (0,24)*{} };
{\ar@{-}_{g \gamma} (25,48)*{}; (50,24)*{} };
(25,24)*{H};
{\ar@{.} (6,0.5)*{}; (6,17.5)*{} };
{\ar@{.} (13,0.5)*{}; (13,11)*{} };
{\ar@{.} (19,0.5)*{}; (19,5)*{} };
{\ar@{.} (31,0.5)*{}; (31,5)*{} };
{\ar@{.} (38,0.5)*{}; (38,11)*{} };
{\ar@{.} (44,0.5)*{}; (44,17.5)*{} };
{\ar@{.} (6,47.5)*{}; (6,30.5)*{} };
{\ar@{.} (13,47.5)*{}; (13,37)*{} };
{\ar@{.} (19,47.5)*{}; (19,43)*{} };
{\ar@{.} (31,47.5)*{}; (31,43)*{} };
{\ar@{.} (38,47.5)*{}; (38,37)*{} };
{\ar@{.} (44,47.5)*{}; (44,30.5)*{} };
\endxy
\]
Both transformation axioms are straightforward to check.  The unit axiom follows from the fact that $H \circ (\gamma \times 1)$ (without the alterations as above) is horizontally constant if $\gamma$ is a constant path.  Thus when $\gamma$ is a constant path, $\Pi_{2}H_{\gamma}$ is homotopic to a composite of unit isomorphisms.  The associativity axiom follows from the fact that the map $H \circ (\delta \gamma \times 1)$, where $\delta \gamma$ here indicates the composite of paths in $\Pi_{2}X$, is the map pictured below.
\[
\xy
{\ar@{-}^{g(\gamma)} (0,0)*+{}; (10,0)*+{} };
{\ar@{-}^{g(\delta)} (10,0)*+{}; (20,0)*+{} };
{\ar@{-}^{H_{y}} (20,0)*+{}; (20,-12)*+{} };
{\ar@{-}_{H_{x}} (0,0)*+{}; (0,-12)*+{} };
{\ar@{-}_{f(\gamma)} (0,-12)*+{}; (10,-12)*+{} };
{\ar@{-}_{f(\delta)} (10,-12)*+{}; (20,-12)*+{} };
(5,-6)*{\scriptstyle H}; (15,-6)*{\scriptstyle H};
{\ar@{-} (10,0)*{}; (10,-12)*{} }
\endxy
\]

Next we define $\Pi_{2}[\alpha]$ for a homotopy class of homotopies $[\alpha]:H \Rrightarrow J$.  Taking the representative $\alpha$ for the class gives a map
\[
\alpha:X \times I \times I \rightarrow Y.
\]
Restricting to the point $x$ gives a continuous map $\alpha_{x}: I \times I \rightarrow Y$ which is a homotopy between $\alpha_{x}(-,0) = H_{x}$ and $\alpha_{x}(-,1) = J_{x}$.  By the definition of the cell $[\alpha]$, we compute that
\[
\alpha_{x}(0,-) = f(x) \qquad \alpha_{x}(1,-) = g(x),
\]
so $\alpha_{x}$ is a 2-cell in $\Pi_{2}Y$.  By definition, it is independent of the choice of representative for $[\alpha]$.  We can now check that the assignment $x \mapsto \alpha_{x}$ satisfies the requirements to be a modification.  This amounts to showing that the following two maps are homotopic, fixing the boundary; here any two-dimensional regions which are unmarked are vertically constant.
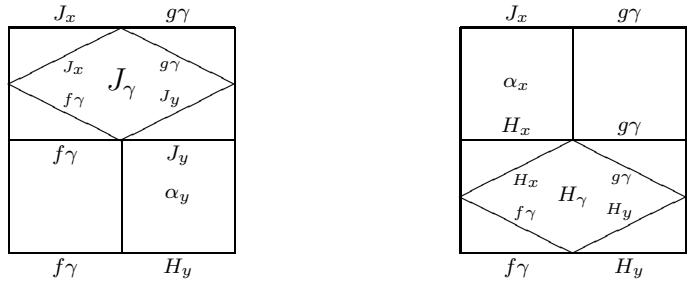
\[
\xy
{\ar@{-}_{\scriptstyle f \gamma} (0,0)*{}; (15,0)*{} };
{\ar@{-}_{\scriptstyle H_{y}} (15,0)*{}; (30,0)*{} };
{\ar@{-}  (0,0)*{}; (0,30)*{} };
{\ar@{-}^{\scriptstyle J_{x}} (0,30)*{}; (15,30)*{} };
{\ar@{-}^{\scriptstyle g \gamma} (15,30)*{}; (30,30)*{} };
{\ar@{-} (30,0)*{}; (30,30)*{} };
{\ar@{-} (15,0)*{}; (15,15)*{} };
{\ar@{-}_{\scriptstyle f \gamma} (0,15)*{}; (15,15)*{} };
{\ar@{-}_{\scriptstyle J_{y}} (15,15)*{}; (30,15)*{} };
{\ar@{-}^{\scriptscriptstyle f \gamma} (0,22.5)*{}; (15,15)*{} };
{\ar@{-}^{\scriptscriptstyle J_{y}} (15,15)*{}; (30,22.5)*{} };
{\ar@{-}_{\scriptscriptstyle J_{x}} (0,22.5)*{}; (15,30)*{} };
{\ar@{-}_{\scriptscriptstyle g \gamma} (15,30)*{}; (30,22.5)*{} };
(22.5, 7.5)*{\scriptstyle \alpha_{y}}; (15,22.5)*{J_{\scriptstyle \gamma}};
{\ar@{-}_{\scriptstyle f \gamma} (60,0)*{}; (75,0)*{} };
{\ar@{-}_{\scriptstyle H_{y}} (75,0)*{}; (90,0)*{} };
{\ar@{-}  (60,0)*{}; (60,30)*{} };
{\ar@{-}^{\scriptstyle J_{x}} (60,30)*{}; (75,30)*{} };
{\ar@{-}^{\scriptstyle g \gamma} (75,30)*{}; (90,30)*{} };
{\ar@{-} (90,0)*{}; (90,30)*{} };
{\ar@{-}^{\scriptscriptstyle f \gamma} (60,7.5)*{}; (75,0)*{} };
{\ar@{-}^{\scriptscriptstyle H_{y}} (75,0)*{}; (90,7.5)*{} };
{\ar@{-}_{\scriptscriptstyle H_{x}} (60,7.5)*{}; (75,15)*{} };
{\ar@{-}_{\scriptscriptstyle g \gamma} (75,15)*{}; (90,7.5)*{} };
{\ar@{-}^{\scriptstyle H_{x}} (60,15)*{}; (75,15)*{} };
{\ar@{-}^{\scriptstyle g \gamma} (75,15)*{}; (90,15)*{} };
{\ar@{-} (75,15)*{}; (75,30)*{} };
(67.5,22.5)*{\scriptstyle \alpha_{x}}; (75,7.5)*{\scriptstyle H_{\gamma}}
\endxy
\]
This follows from the fact that the maps below are homotopic fixing the top and bottom boundaries since taken together these maps form four of the six faces of the image of the cube $\alpha \circ (\gamma \times 1 \times 1):I^{3} \rightarrow Y$.
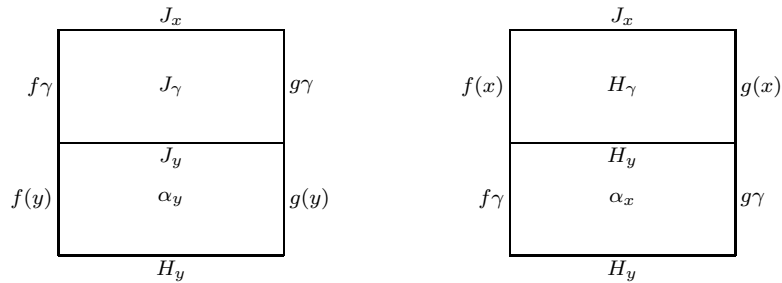
\[
\xy
{\ar@{-}_{\scriptstyle H_{y}} (0,0)*{}; (30,0)*{} };
{\ar@{-}_{\scriptstyle g(y)} (30,0)*{}; (30,15)*{} };
{\ar@{-}_{\scriptstyle g \gamma} (30,15)*{}; (30,30)*{} };
{\ar@{-}^{\scriptstyle f(y)} (0,0)*{}; (0,15)*{} };
{\ar@{-}^{\scriptstyle f \gamma} (0,15)*{}; (0,30)*{} };
{\ar@{-}^{\scriptstyle J_{x}} (0,30)*{}; (30,30)*{} };
{\ar@{-}_{\scriptstyle J_{y}} (0,15)*{}; (30,15)*{} };
(15,7.5)*{\scriptstyle \alpha_{y}}; (15,22.5)*{\scriptstyle J_{\gamma}};
{\ar@{-}_{\scriptstyle H_{y}} (60,0)*{}; (90,0)*{} };
{\ar@{-}_{\scriptstyle g \gamma} (90,0)*{}; (90,15)*{} };
{\ar@{-}_{\scriptstyle g(x)} (90,15)*{}; (90,30)*{} };
{\ar@{-}^{\scriptstyle f \gamma} (60,0)*{}; (60,15)*{} };
{\ar@{-}^{\scriptstyle f(x)} (60,15)*{}; (60,30)*{} };
{\ar@{-}^{\scriptstyle J_{x}} (60,30)*{}; (90,30)*{} };
{\ar@{-}_{\scriptstyle H_{y}} (60,15)*{}; (90,15)*{} };
(75,7.5)*{\scriptstyle \alpha_{x}}; (75,22.5)*{\scriptstyle H_{\gamma}};
\endxy
\]
This completes the description of $\Pi_{2}$ on the cells of $\mb{Top}$, so we now are in a position to prove that $\Pi_{2}$ is a functor.

\begin{thm}
The map on underlying 3-globular sets given above can be given the structure of a functor of tricategories
\[
\Pi_{2}: \mb{Top} \rightarrow \mb{Bicat}.
\]
\end{thm}

$\mb{Proof.}$
To give this map of underlying 3-globular sets the structure of a strict functor between tricategories, we need only check that this map coherently preserves all units and compositions.  We begin with 1-cells and work our way up.

By construction, $\Pi_{2}f$ is a strict functor.  Thus given a composable pair $(g,f)$, both $\Pi_{2}(gf)$ and $\Pi_{2}g \circ \Pi_{2}f$ are strict functors which agree on cells, so $\Pi_{2}(gf) = \Pi_{2}g \circ \Pi_{2}f$.  Additionally, $\Pi_{2}1$ is the strict functor which is the identity on cells, so it is the identity functor.  Thus $\Pi_{2}$ strictly preserves 1-cell composition and units.

We now show that the map $\Pi_{2}: \mb{Top}(X,Y) \rightarrow \mb{Bicat}(\Pi_{2}X, \Pi_{2}Y)$ can be given the structure of a functor between bicategories.  First, it is clear that $\Pi_{2}$ sends an identity 3-cell $1: H \Rrightarrow H$ of $\mb{Top}$ to the identity modification $1_{\Pi_{2}H}$.  Second, we have that $\Pi_{2}[\beta] \circ \Pi_{2}[\alpha] = \Pi_{2}[\beta \alpha]$ since both sides are obtained by composing the 2-cells $\beta_{x}$ and $\alpha_{x}$ vertically in $\Pi_{2}Y$. Now let $H:f \Rightarrow f'$ and $J:f' \Rightarrow f''$ be a pair of composable 1-cells in $\mb{Top}(X,Y)$.  The definitions immediately imply that $\Pi_{2}(JH)_{x} = \Pi_{2}J_{x} \circ \Pi_{2}H_{x}$, and that the maps $\Pi_{2}(JH)_{\gamma}$ and $(\Pi_{2}J \circ \Pi_{2}H)_{\gamma}$ are homotopic by a homotopy fixing the boundary square so $[\Pi_{2}(JH)_{\gamma}] = [(\Pi_{2}J \circ \Pi_{2}H)_{\gamma}]$ and therefore $\Pi_{2}(JH) = \Pi_{2}J \circ \Pi_{2}H$.  Finally, consider the identity homotopy $1_{f}:f \Rightarrow f$.  The transformation $\Pi_{2}1_{f}$ has its component at $x$ the constant path, and the map $(\Pi_{2}1_{f})_{\gamma}$ is easily seen to be homotopic, fixing the boundary, to the composite
\[
1_{\gamma(1)} f(\gamma) \stackrel{l}{\longrightarrow} f(\gamma) \stackrel{r^{-1}}{\longrightarrow} f(\gamma) 1_{\gamma(0)},
\]
so $\Pi_{2}$ sends the identity $1_{f}$ to the identity transformation $\Pi_{2}f \Rightarrow \Pi_{2}f$.  This shows that the map $\Pi_{2}: \mb{Top}(X,Y) \rightarrow \mb{Bicat}(\Pi_{2}X, \Pi_{2}Y)$ is a strict functor of bicategories.

Next, we give the rest of the data for a functor between tricategories, postponing any axioms until afterwards.  We first complete the definition of the adjoint equivalence
\[
\chi:\Pi_{2} \circ \otimes_{\mb{Top}} \simeq \otimes_{\mb{Bicat}} \circ (\Pi_{2} \times \Pi_{2}).
\]
We know that these two functors agree on objects, so we define the component of $\chi$ at $(g,f)$ to be the identity.  Now $\Pi_{2}(J \otimes H)_{x}$ is given by the path $J(H(x,s),s)$, while $(\Pi_{2}J \otimes \Pi_{2}H)_{x}$ is given by the path below.
\[
(\Pi_{2}J \otimes \Pi_{2}H)_{x}(s) = \left\{ \begin{array}{cl}
J(H(x,0),2s) & \quad 0 \leq s \leq \frac{1}{2} \\
J(H(x,2s-1),1) & \quad \frac{1}{2} \leq s \leq 1 \end{array} \right.
\]
The transformation $\chi$ therefore has a component at $(J,H)$ of the form
\[
1 \circ (\Pi_{2}J \otimes \Pi_{2}H) \stackrel{\chi}{\Rrightarrow} \Pi_{2}(J \otimes H) \circ 1.
\]
This is uniquely determined by 2-cells $(\Pi_{2}J \otimes \Pi_{2}H)_{x} \Rightarrow \Pi_{2}(J \otimes H)_{x}$ in $\Pi_{2}Y$ which we define by the formula below.
\[
\chi_{(J,H)}(x,s,t) = \left\{ \begin{array}{cl}
J(H(x,st), 2s-st) & \quad 0 \leq s \leq \frac{1}{2} \\
J(H(x,2s-st+t-1), st-t+1) & \quad \frac{1}{2} \leq s \leq 1 \end{array} \right.
\]

The data for a functor also includes an adjoint equivalence $\bs{\iota}$ between the unit in the target and the image of the unit in the source.  It is easy to check that $\Pi_{2}$ sends the identity map $1:X \rightarrow X$ in $\mb{Top}$ to the identity functor $1: \Pi_{2}X \rightarrow \Pi_{2}X$ in $\mb{Bicat}$, so $\mb{\iota}$ is defined to be the identity adjoint equivalence.

Finally, there are invertible modifications $\omega, \delta, \gamma$.  Each of these has source and target a composite of coherence 1-cells, so we define them in each case to be given by unique coherence 2-cells.  We must check that this collection of data satisfies the axioms to be a modification, but in each case this follows using simple reparametrization homotopies.  In addition, these definitions immediately imply the two functor axioms.

The only thing left to check is that the data given for $\chi$ actually produces a transformation.  First, we must check that this definition makes $x \mapsto \chi_{(J,H)}(x,s,t)$ a modification.  Given a path $\gamma: x \rightarrow y$ in $X$, we are required to check that two different composites of 2-cells in $\Pi_{2}Y$ are equal. This is accomplished by taking the homotopy $\Gamma$ below and modifying it as in the construction of $\Pi_{2}$.
\[
\Gamma(s,t,r) = \left\{ \begin{array}{cl}
J(H(y,2st),2s-2st) & 0 \leq s \leq \frac{1}{2}, 0 \leq t \leq \frac{r}{2} \\
J(H(\gamma(1+r-2t), rs), 2s-rs) & 0 \leq s \leq \frac{1}{2}, \frac{r}{2} \leq t \leq \frac{r+1}{2} \\
J(H(x, 2st-s), 3s-2st) & 0 \leq s \leq \frac{1}{2}, \frac{r+1}{2} \leq t \leq 1 \\
J(H(y, 2s-2st+2t-1), 2st-2t+1) & \frac{1}{2} \leq s \leq 1, 0 \leq t \leq \frac{r}{2} \\
J(H(\gamma(1+r-2t), 2s-rs+r-1), rs-r+1) & \frac{1}{2} \leq s \leq 1, \frac{r}{2} \leq t \leq \frac{r+1}{2} \\
J(H(x,3s-2st+2t-2), 2st-s-2t+2) & \frac{1}{2} \leq s \leq 1, \frac{r+1}{2} \leq t \leq 1 \end{array} \right.
\]
Then we must show that the 2-cell $\chi_{(1,1)}(x,s,t)$ gives the identity and that the horizontal composite $\chi_{(J',H')} * \chi_{(J,H)}$ equals $\chi_{(J'J, H'H)}$ precomposed with a naturality 2-cell for $\Pi_{2}J'$ with respect to $(\Pi_{2}H)_{x}$.  The first of these is completely trivial and the second can be verified using a tedious but straightforward contracting homotopy.  $\Box$

\subsection{Monoidal fundamental 2-groupoids}

We are now in a position to prove the two main results of this section, namely that taking the fundamental 2-groupoids of algebras for the operads $\mathcal{C}_{1}$ and $\mathcal{C}_{2}$ yield monoidal and braided monoidal bicategories.  While the braided monoidal case is the one of greater interest, we use the plain monoidal case to illustrate the main ideas.

The category of topological spaces has Cartesian products, as does the category of bicategories and weak functors.  It is clear that the functor $\Pi_{2}$ between tricategories defined above restricts to a functor between ordinary categories $\Pi_{2}:Top \rightarrow Bicat$.

\begin{lem}
The functor (of ordinary categories) $\Pi_{2}:Top \rightarrow Bicat$ is monoidal with respect to the Cartesian structures.
\end{lem}

The key construction we require in this section is given in the following simple lemma.

\begin{lem}\label{operadpts}
Let $\mathcal{P}$ be an operad in the category of topological spaces, and let $X$ be an algebra.  Then every $p \in \mathcal{P}(n)$ gives a map $\mu_{p}:X^{n} \rightarrow X$, and every path $\gamma:I \rightarrow \mathcal{P}(n)$ gives a homotopy
\[
\tilde{\gamma}: \mu_{\gamma(0)} \Rightarrow \mu_{\gamma(1)}.
\]
\end{lem}

\begin{thm}
Let $X$ be an algebra for the operad $\mathcal{C}_{1}$.  Then $\Pi_{2}X$ has the structure of a monoidal bicategory.
\end{thm}
$\mb{Proof.}$
The underlying bicategory of $\Pi_{2}X$ is the bicategory of the same name constructed above.  The tensor product
\[
\otimes: \Pi_{2}X \times \Pi_{2}X \rightarrow \Pi_{2}X
\]
is given by the functor
\[
\Pi_{2}X \times \Pi_{2}X \cong \Pi_{2}(X \times X) \stackrel{\Pi_{2}\mu_{m}}{\longrightarrow} \Pi_{2}X
\]
where $m$ is the element of $\mathcal{C}_{1}(2)$ given by the pair of 1-cubes (i.e., intervals) $(\frac{1}{5}, \frac{2}{5})$ and $(\frac{3}{5}, \frac{4}{5})$, and $\mu_{m}$ is the map $X^{2} \rightarrow X$ given by Lemma \ref{operadpts}.  The unit functor is
\[
\Pi_{2}\mu_{i}: \Pi_{2} \{ * \} \rightarrow \Pi_{2}X,
\]
where $i$ is the unique point in $\mathcal{C}_{1}(0)$.

For the associativity adjoint equivalence, first note that $(x \otimes y) \otimes z$ is the point in $X$ given by the formula
\[
\mu_{m}(\mu_{m}(x, y), z)
\]
while $x \otimes (y \otimes z)$ is the point given by the formula
\[
\mu_{m}(x, \mu_{m}(y, z)).
\]
We also have the operad multiplication, giving maps
\[
\begin{array}{c}
\sigma: \mathcal{C}_{1}(2) \times \mathcal{C}_{1}(2) \times \mathcal{C}_{1}(1) \rightarrow \mathcal{C}_{1}(3) \\
\sigma': \mathcal{C}_{1}(2) \times \mathcal{C}_{1}(1) \times \mathcal{C}_{1}(2) \rightarrow \mathcal{C}_{1}(3).
\end{array}
\]
Thus we see that $(x \otimes y) \otimes z$ is given by evaluating $\mathcal{C}_{1}(3) \times X^{3} \rightarrow X$ at $\big( \sigma(m,m,1), x, y, z \big)$, and  $x \otimes (y \otimes z)$ is given by evaluating the same map $\mathcal{C}_{1}(3) \times X^{3} \rightarrow X$ at $\big( \sigma'(m,1,m), x, y, z \big)$.  Therefore, to give the associativity adjoint equivalence it suffices, by Lemma \ref{operadpts}, to give a path in $\mathcal{C}_{1}(3)$ from $\sigma(m,m,1)$ to $\sigma'(m,1,m)$.  Since a little 1-cube is determined by its center and its length, we provide a path using that information writing this as a triple of points on the real line.  The path, given by the map $\alpha:I \rightarrow \mathcal{C}_{1}(3)$, is defined to have its center at $(\frac{13+2t}{50}, \frac{17+16t}{50}, \frac{35+2t}{50})$ with lengths $(\frac{1+4t}{25}, \frac{1}{25}, \frac{5-4t}{25})$.



The component $a_{xyz}$ is then the 1-cell represented by the path $\alpha:I \rightarrow X$ which is the composite
\[
I \stackrel{1 \times (x,y,z)}{\longrightarrow} I \times X^{3} \stackrel{\alpha \times 1}{\longrightarrow} \mathcal{C}_{1}(3) \times X^{3} \longrightarrow X.
\]
Now given 1-cells in $\Pi_{2}X$ represented by paths $f:x \rightarrow x', g: y \rightarrow y', h: z \rightarrow z'$, we must construct a 2-cell $a_{fgh}$ in $\Pi_{2}X$.  This 2-cell is obtained just as we did above for the 2-cell isomorphism data for transformations by taking the map
\[
I \times I \stackrel{\alpha \times (f \times g \times h)}{\longrightarrow} \mathcal{C}_{1}(3) \times X^{3} \longrightarrow X
\]
and modifying it to be a 2-cell as in the construction of $\Pi_{2}X$ in Section 3.1.  The pseudoinverse $a^{\centerdot}$ is then given by the same procedure using the inverse path, $\alpha^{-1}$, and the unit and counit of this adjoint equivalence are the canonical contracting homotopies $\alpha \circ \alpha^{-1} \simeq 1, \alpha^{-1} \circ \alpha \simeq 1$.  An analogous construction gives the unit adjoint equivalences.

For the four invertible modifications, note that each $\mathcal{C}_{1}(j)$ is a disjoint union of contractible spaces.  Thus there is a unique homotopy class of homotopies between the source and target paths.  Additionally, this means that the monoidal bicategory axioms follow trivially.
$\Box$

Now we give the construction of, for an algebra $X$ over the little 2-cubes operad, a braided monoidal bicategory $\Pi_{2}X$.  The reader should note that the construction here does not precisely extend the one given for $\mathcal{C}_{1}$-algebras in the following sense.  There is a map of operads $\mathcal{C}_{1} \rightarrow \mathcal{C}_{2}$ which sends a collection of little 1-cubes $(\gamma_{i}) \in \mathcal{C}_{1}(n)$ to the collection $(\gamma_{i} \times J) \in \mathcal{C}_{2}(n)$ where $J$ is the open unit interval.  This gives every $\mathcal{C}_{2}$-algebra $X$ the structure of a $\mathcal{C}_{1}$-algebra by restriction, which we will denote $RX$.  From the proof below, it is immediate that the underlying monoidal bicategory of $\Pi_{2}X$ is not equal to $\Pi_{2}RX$, but only monoidally biequivalent to it by a functor which is the identity on underlying cells.  On the other hand, the construction is largely identical to the one for the monoidal structure on a $\mathcal{C}_{1}$-algebra, so the bulk of the proof focuses on the braiding itself.

\begin{thm}
Let $X$ be an algebra for the operad $\mathcal{C}_{2}$.  Then $\Pi_{2}X$ has the structure of a braided monoidal bicategory.
\end{thm}
$\mb{Proof.}$
The tensor for $\Pi_{2}X$ is now given by restricting the operad action to the element $m \in \mathcal{C}_{2}(2)$ given by the pair of little 2-cubes
\[
\Big( (\frac{1}{5}, \frac{2}{5}) \times (\frac{2}{5}, \frac{3}{5}), (\frac{3}{5}, \frac{4}{5}) \times (\frac{2}{5}, \frac{3}{5}) \Big)
\]
using Lemma \ref{operadpts}.  It will be useful to describe little 2-cubes, and in particular paths in the space $\mathcal{C}_{2}(k)$, by stating the center and size of each little 2-cube.  Using this definition of tensor product, we compute that the triple tensor $(xy)z$ is given by restricting the action of operad to the point in $\mathcal{C}_{2}(3)$ with centers at
\[
\Big( (\frac{13}{50}, \frac{1}{2}), (\frac{17}{50}, \frac{1}{2}), (\frac{35}{50}, \frac{1}{2}) \Big)
\]
and side lengths $(\frac{1}{25}, \frac{1}{25}, \frac{1}{5})$, while the triple tensor $x(yz)$ is given by restricting to the point with centers at
\[
\Big( (\frac{15}{50}, \frac{1}{2}), (\frac{33}{50}, \frac{1}{2}), (\frac{37}{50}, \frac{1}{2}) \Big)
\]
 and side lengths $(\frac{1}{5}, \frac{1}{25}, \frac{1}{25})$.  The associativity equivalence is given by the path with centers
\[
\Big( (\frac{13+2t}{50}, \frac{1}{2}), (\frac{17+16t}{50}, \frac{1}{2}), (\frac{35+2t}{50}, \frac{1}{2}) \Big)
\]
and side lengths $(\frac{1+4t}{25}, \frac{1}{25}, \frac{5-4t}{25})$.  The rest of the adjoint equivalence is then constructed using the reverse homotopy and the obvious contracting maps, and the left and right unit adjoint equivalences are given similarly.

The braid is given by using using the path from $x \otimes y$ to $y \otimes x$ that we describe now.  At time $t$, the braid is the pair of little 2-cubes $(b_{1}(t), b_{2}(t))$ where each $b_{i}(t)$ is size $\frac{1}{5} \times \frac{1}{5}$ with $b_{1}(t)$ centered at
\[
\Big( \frac{1}{2} + \frac{1}{5}\cos (\pi + \pi t), \frac{1}{2} + \frac{1}{5}\sin (\pi + \pi t) \Big)
\]
and $b_{2}(t)$ centered at
\[
\Big(\frac{1}{2} + \frac{1}{5}\cos (\pi t), \frac{1}{2} + \frac{1}{5}\sin (\pi t) \Big).
\]
The rest of the adjoint equivalence $\mb{R}$ is given using the reverse path and the obvious contracting homotopies.

The final data required is that of two modifications $R_{(-|-,-)}, R_{(-,-|-)}$.  We will give explicit formulas for the first of these, the second is constructed in precisely the same fashion.  We are required to give components $R_{(x|y,z)}$, and to do that we provide a map $D^{2} \rightarrow \mathcal{C}_{2}(3)$; this map will give a homotopy between the source and target paths in $\mathcal{C}_{2}(3)$, so by the same argument as in Lemma \ref{operadpts} will produce a 2-cell in $\Pi_{2}X$.

Both source and target paths have three little 2-cubes, which we give now.  Since both the source and target 1-cells of $R_{(x|y,z)}$ are the composite of three generating 1-cells, we must compose three different paths in $\Pi_{2}X$.  We ignore this detail here as it is irrelevant, and instead replace these composed paths each with a path of length three in which each of the generating 1-cell paths is traversed in one unit of time.  We also ignore the side lengths of the cubes, as these can always be made small enough to be irrelevant, so we denote paths only by where the centers of each cube are located.  Using these conventions, the source is the collection of the three paths given below.
\[
\gamma_{1}(t) = \left\{ \begin{array}{cl}
\Big( \frac{3}{10} + \frac{1}{25}\cos (\pi + \pi t), \frac{1}{2} + \frac{1}{25} \sin (\pi + \pi t) \Big)& 0 \leq t \leq 1 \\
\Big( \frac{17 + 16(2t-1)}{50}, \frac{1}{2} \Big) & 1 \leq t \leq 2 \\
\Big( \frac{7}{10} + \frac{1}{25}\cos (\pi + \pi (3t-2)), \frac{1}{2} + \frac{1}{25} \sin ( \pi + \pi(3t-2)) \Big) & 2 \leq t \leq 3 \end{array} \right.
\]
\[
\gamma_{2}(t) = \left\{ \begin{array}{cl}
\Big( \frac{3}{10} + \frac{1}{25}\cos (\pi t), \frac{1}{2} + \frac{1}{25} \sin (\pi t)\Big) & 0 \leq t \leq 1 \\
\Big( \frac{13+2t}{50}, \frac{1}{2} \Big) & 1 \leq t \leq 2 \\
\Big( \frac{3}{10}, \frac{1}{2} \Big) & 2 \leq t \leq 3 \end{array} \right.
\]
\[
\gamma_{3}(t) = \left\{ \begin{array}{cl}
\Big( \frac{7}{10}, \frac{1}{2} \Big) & 0 \leq t \leq 1 \\
\Big( \frac{35+2(2t-1)}{50}, \frac{1}{2} \Big) & 1 \leq t \leq 2 \\
\Big( \frac{7}{10} + \frac{1}{25} \cos (\pi (3t-2)), \frac{1}{2} + \frac{1}{25} \sin (\pi (3t-2)) \Big) & 2 \leq t \leq 3 \end{array} \right.
\]
We can similarly compute that the target is the collection of the three paths given here.
\[
\gamma_{1}'(t) = \left\{ \begin{array}{cl}
\Big( \frac{13+2t}{50}, \frac{1}{2} \Big) & 0 \leq t \leq 1 \\
\Big( \frac{1}{2} + \frac{1}{5} \cos (\pi + \pi (2t-1)), \frac{1}{2} + \frac{1}{5} \sin ( \pi + \pi (2t-1)) \Big) & 1 \leq t \leq 2 \\
\Big( \frac{35 + 2(3t-2)}{50}, \frac{1}{2} \Big) & 2 \leq t \leq 3 \end{array} \right.
\]
\[
\gamma_{2}'(t) = \left\{ \begin{array}{cl}
\Big( \frac{17+16t}{50}, \frac{1}{2} \Big) & 0 \leq t \leq 1 \\
\Big( \frac{23}{50} + \frac{1}{5} \cos ( \pi (2t-1)), \frac{1}{2} + \frac{1}{5} \sin ( \pi (2t-1)) \Big) & 1 \leq t \leq 2 \\
\Big( \frac{13+2(3t-2)}{50}, \frac{1}{2} \Big) & 2 \leq t \leq 3 \end{array} \right.
\]
\[
\gamma_{3}'(t) = \left\{ \begin{array}{cl}
\Big( \frac{35+2t}{50}, \frac{1}{2} \Big) & 0 \leq t \leq 1 \\
\Big( \frac{27}{50} + \frac{1}{5} \cos ( \pi (2t-1)), \frac{1}{2} + \frac{1}{5} \sin ( \pi (2t-1)) \Big) & 1 \leq t \leq 2 \\
\Big( \frac{17+16(3t-2)}{50}, \frac{1}{2} \Big) & 2 \leq t \leq 3 \end{array} \right.
\]
Both of these are homotopic to the collection of three paths $\delta$ given below (once again of length three) by the obvious linear homotopies, and composing the homotopies $\gamma \simeq \delta \simeq \gamma'$ gives the required map.
\[
\begin{array}{rcl}
\delta_{1}(t) & = & \Big( \frac{1}{2} + \frac{6}{25} \cos (\pi + \frac{\pi}{3}t), \frac{1}{2} + \frac{6}{25} \sin (\pi + \frac{\pi}{3}t) \Big) \\
\delta_{2}(t) & = & \Big( \frac{3}{10} + \frac{2}{25} \cos(\frac{\pi}{3}t), \frac{1}{2} + \frac{2}{25} \sin(\frac{\pi}{3}t) \Big) \\
\delta_{3}(t) & = & \Big( \frac{17}{25} + \frac{1}{50}\cos(\frac{\pi}{3}t), \frac{1}{2} + \frac{1}{50}\sin(\frac{\pi}{3}t) \Big)
\end{array}
\]

There are now axioms to check to show that this collection of data gives a braided monoidal bicategory.  In each axiom, the 2-cell pastings to be shown equal are given by maps $D^{2} \rightarrow \mathcal{C}_{2}(k)$ which are then used along with the operad action to define the actual pasting.  Since every $\mathcal{C}_{2}(k)$ has trivial homotopy groups above dimension one, every such pair of maps $D^{2} \rightarrow \mathcal{C}_{2}(k)$ with the same boundary, such as those arising from the braided monoidal bicategory axioms, are necessarily homotopic in $\mathcal{C}_{2}(k)$.  These homotopies then show that the two pastings required to be equal for an axiom to hold are in fact equal, as equality of 2-cells in $\Pi_{2}X$ is exactly given by such a homotopy. $\Box$

\subsection{Lifting structures and homotopy invariance}

The final tool needed for coherence is a homotopy-invariance result.  Using this theorem we will be able to give certain fundamental 2-groupoids extra algebraic structure in the next section.  We note that it would be possible to omit this discussion if we had chosen a cofibrant braided operad (see \cite{fie} for a discussion of braided operads) instead of the little 2-cubes operad, but we chose to retain the original operad since using the little 2-cubes operad made the construction of the braiding on the fundamental 2-groupoid of an algebra transparent.  Furthermore, our proof of coherence explicitly uses the relationship between surface braids and configuration spaces, making it advantageous to choose an operad whose spaces can be easily compared with configuration spaces.  One preliminary result is needed first.

\begin{thm}[Transfer of structure]
Let $X$ be an object in $\mb{MonBicat}$ (resp., $\mb{BrMonBicat}$), and let $Y$ be any bicategory.  Let $f:Y \rightarrow X$ be a biequivalence from $Y$ to the underlying bicategory of $X$.  Then $Y$ can be given a monoidal (resp., braided monoidal) structure such that $f$ is a monoidal (resp., braided monoidal) functor.  In fact, $f$ can be completed to give a biadjoint biequivalence between $Y$ and $X$ in $\mb{MonBicat}$ (resp., $\mb{BrMonBicat}$).
\end{thm}
$\mb{Proof.}$
This is Theorem 5.1 and Remark 5.3 of \cite{g2}.  $\Box$

\begin{cor}[Homotopy invariance of structure]
Let $X$ be an algebra for the little $n$-cubes operad, $n=1,2$, and let $f:Y \rightarrow X$ be a homotopy equivalence.  Then $\Pi_{2}Y$ can be given the structure of an object in $\mb{MonBicat}$ for $n=1$ or $\mb{BrMonBicat}$ for $n=2$ such that $\Pi_{2}f$ is a monoidal (when $n=1$) or braided monoidal (when $n=2$) biequivalence.
\end{cor}
$\mb{Proof.}$
First, note that a 1-cell $f$ in $\mb{Top}$ is a biequivalence, and hence part of an internal biadjoint biequivalence by \cite{g2}, if and only if it is a homotopy equivalence.  Since every functor sends biadjoint biequivalences to biadjoint biequivalences,
\[
\Pi_{2}f: \Pi_{2}Y \rightarrow \Pi_{2}X
\]
will then be a biadjoint biequivalence.  Using the previous theorem, we can then lift the braided monoidal structure from $\Pi_{2}X$ via $\Pi_{2}f$.  $\Box$

\section{Free monoidal bicategories}

This section shows how free monoidal and braided monoidal bicategories can be interpreted topologically.  In particular, we show that the free $n$-tuply monoidal bicategory on one object is $n$-tuply monoidally biequivalent to the fundamental 2-groupoid $\Pi_{2} \big( \coprod_{k} B(k, \mathbb{R}^{n}) \big)$ for the cases $n=1$ and $n=2$.  We first review all of the free constructions required.  Next we briefly discuss the case $n=1$.  This case concerns monoidal bicategories for which a coherence theorem is already known, but we use it to derive a topological interpretation of free monoidal bicategories in order to state the kind of theorem that we call coherence for braided monoidal bicategories.  Finally, we prove the main result for the case $n=2$.

\subsection{Free structures}

Here we will review the construction of free monoidal bicategories, free $\textbf{Gray}$-monoids, and free braided monoidal bicategories.  The first two of these stuctures are studied in \cite{g1}, but the third is new.  All of these objects are constructed in a similar fashion:  first we inductively construct all of the required cells, and then identify cells as required by the necessary axioms.  The free constructions we present here are all left adjoints to the forgetful functor to some category of underlying data; in each case, we explain the universal property.

Let $X$ be a 2-category.  The free $\textbf{Gray}$-monoid on $X$, $\mathcal{F}_{Gr}X$, has objects consisting of all finite strings of elements of $X$ including the empty string.  The set of morphisms from one string $\underline{x}$ to another string $\underline{y}$ in $\mathcal{F}_{Gr}X$ is empty if the length of $\underline{x}$ is different from the length of $\underline{y}$, and generated under composition by morphisms of the form
\[
1 \otimes 1 \otimes \cdots \otimes f_{i} \otimes \cdots \otimes 1
\]
when the lengths are the same, and $f_{i}:x_{i} \rightarrow y_{i}$ is a 1-cell of $X$.  The two cells of $\mathcal{F}_{Gr}X$ are generated by the 2-cells of $X$ and new isomorphisms
\[
(f \otimes 1) \circ (1 \otimes g) \cong (1 \otimes g) \circ (f \otimes 1),
\]
subject to the usual 2-category axioms along with the new $\mb{Gray}$-monoid axioms.  The reader not familiar with these axioms is invited to consult \cite{g1} or \cite{gps}.  It should be noted that if $X$ is a set seen as a 2-category with only identity 1- and 2-cells, then $\mathcal{F}_{Gr}X$ is the free monoid on $X$ treated as a discrete monoidal 2-category.

Let $\mb{2Cat}$ denote the category of 2-categories and 2-functors between. Let $\mb{GrayMon}$ denote the category of $\mb{Gray}$-monoids and $\mb{Gray}$-functors between them; these are the 2-functors between the underlying 2-categories which strictly preserve the multiplication and unit.  We have an obvious forgetful functor $\mb{GrayMon} \rightarrow \mb{2Cat}$ which sends each $\mb{Gray}$-monoid to its underlying 2-category.  The universal property of the free $\mb{Gray}$-monoid construction is then expressed by the following proposition.

\begin{prop}
The free $\mb{Gray}$-monoid functor $\mathcal{F}_{Gr}: \mb{2Cat} \rightarrow \mb{GrayMon}$ is left adjoint to the forgetful functor.
\end{prop}

The free monoidal bicategory is constructed in much the same way as the free $\mb{Gray}$-monoid, but has more generating cells.  Let $X$ be a bicategory.  The objects of $\mathcal{F}X$ are generated by the elements of $X$ and a new unit object $I$ by taking binary tensor products; this gives the tensor product $\otimes$ on $\mathcal{F}X$ which we will often omit.  The 1-cells of $\mathcal{F}X$ are generated under tensor and composition by the 1-cells of $X$ and new 1-cells
\[
\begin{array}{ll}
a_{xyz}:(xy)z \rightarrow x(yz) & \quad a_{xyz}^{\centerdot}:x(yz) \rightarrow (xy)z \\
l_{x}:Ix \rightarrow x & \quad l_{x}^{\centerdot}:x \rightarrow Ix \\
r_{x}:xI \rightarrow x & \quad r_{x}^{\centerdot}:x \rightarrow xI
\end{array}
\]
for all objects $x,y,z$.  The 2-cells of $\mathcal{F}X$ are generated by the 2-cells of $X$, isomorphisms
\[
\begin{array}{ll}
\eta_{a}:1_{(xy)z} \Rightarrow a^{\centerdot}\circ a & \quad \varepsilon_{a}:a\circ a^{\centerdot} \Rightarrow 1_{x(yz)} \\
\eta_{l}:1_{Ix} \Rightarrow l^{\centerdot}\circ l & \quad \varepsilon_{l}:l\circ l^{\centerdot} \Rightarrow 1_{x} \\
\eta_{r}:1_{xI} \Rightarrow r^{\centerdot}\circ r & \quad \varepsilon_{r}:r\circ r^{\centerdot} \Rightarrow 1_{x}
\end{array}
\]
for all objects $x,y,z$, naturality isomorphisms making $a, a^{\centerdot}, l, l^{\centerdot}, r, r^{\centerdot}$ pseudonatural transformations, functoriality isomorphisms making the assignments $(x,y) \mapsto xy$, $* \mapsto I$ weak functors, and four new isomorphisms listed below.
\[
\begin{array}{c}
\pi: (1 \otimes a) \circ a \circ (a \otimes 1) \Rightarrow a \circ a \\
\mu: (1 \otimes l) \circ a \circ r^{\centerdot} \Rightarrow 1 \\
\lambda: l \otimes 1 \Rightarrow l \circ a \\
\rho: 1 \otimes r^{\centerdot} \Rightarrow a \circ r^{\centerdot}
\end{array}
\]
These 2-cells are required to make each of the quadruples
\[
(a, a^{\centerdot}, \eta_{a}, \varepsilon_{a}),  (l, l^{\centerdot}, \eta_{l}, \varepsilon_{l}), (r, r^{\centerdot}, \eta_{r}, \varepsilon_{r})
 \]
into an adjoint equivalence, to satisfy the monoidal bicategory axioms for $\pi, \mu, \lambda, \rho$, and to satisfy the naturality and functoriality axioms mentioned above, so we quotient by the equivalence relation generated by these requirements.

Let $\mb{Bicat}_{s}$ denote the category of bicategories and strict functors between them.  Let $\mb{MonBicat}_{s}$ denote the category of monoidal bicategories and strict monoidal functors between them; these functors have underlying functors of bicategories which are strict, and also strictly preserve all of the monoidal structure.  There is an obvious forgetful functor $\mb{MonBicat}_{s} \rightarrow \mb{Bicat}_{s}$ which forgets the monoidal structure.  The universal property of the free monoidal bicategory functor is expressed by the following proposition.

\begin{prop}
The free monoidal bicategory functor $\mathcal{F}: \mb{Bicat}_{s} \rightarrow \mb{MonBicat}_{s}$ is left adjoint to the forgetful functor.
\end{prop}

The free braided monoidal bicategory on $X$, $\mathcal{F}_{br}X$, is constructed analogously to the free monoidal bicategory but with extra generating 1- and 2-cells and new axioms.  The additional generating 1-cells are the braiding and its pseudoinverse $R_{xy}:xy \rightarrow yx, R_{xy}^{\centerdot}:yx \rightarrow xy$, and the new generating 2-cells are those listed below, together with naturality 2-cells for $R$ and $R^{\centerdot}$ which we do not list.
\[
\begin{array}{ll}
\eta_{R}:1_{xy} \Rightarrow R^{\centerdot}\circ R & \quad \varepsilon_{R}:R\circ R^{\centerdot} \Rightarrow 1_{yx} \\
R_{(x,y|z)}:(R1) \circ a^{\centerdot} \circ (1R) \Rightarrow a^{\centerdot} \circ R \circ a^{\centerdot} & \quad R_{(x|y,z)}: (1R) \circ a \circ (R1) \Rightarrow a \circ R \circ a
\end{array}
\]
Once again, we require $(R, R^{\centerdot}, \eta_{R}, \varepsilon_{R})$ to constitute an adjoint equivalence, that $R$ and $R^{\centerdot}$ are pseudonatural in both variables, and that the braided monoidal bicategory axioms hold for $R_{(-,-|-)}$ and $R_{(-|-,-)}$, all in addition to requiring that the free braided monoidal bicategory is also a monoidal bicategory.

As before, let $\mb{Bicat}_{s}$ denote the category of bicategories and strict functors between them.  Let $\mb{BrMonBicat}_{s}$ denote the category of braided monoidal bicategories and strict braided monoidal functors between them; these functors have underlying functors of bicategories which are strict, and also strictly preserve all of the braided monoidal structure.  There is an obvious forgetful functor $\mb{BrMonBicat}_{s} \rightarrow \mb{Bicat}_{s}$ which forgets the entire braided monoidal structure.  The universal property of the free braided monoidal bicategory functor is expressed by the following proposition.

\begin{prop}
The free braided monoidal bicategory functor $\mathcal{F}_{br}: \mb{Bicat}_{s} \rightarrow \mb{BrMonBicat}_{s}$ is left adjoint to the forgetful functor.
\end{prop}

It should be clear that the forgetful functor $\mb{BrMonBicat}_{s} \rightarrow \mb{Bicat}_{s}$ factors through $\mb{MonBicat}_{s}$.  We could in fact show that forgetting the braiding yields a functor
\[
\mb{BrMonBicat}_{s} \rightarrow \mb{MonBicat}_{s}
\]
which itself has a left adjoint.  This left adjoint takes a monoidal bicategory and freely adds just the braiding to it.  Freely adding a braiding to an already-monoidal bicategory involves adjoining the 1- and 2-cells we have listed above, and then imposing axioms at the level of 2-cells.  These axioms are the braided monoidal bicategory axioms, plus naturality axioms with respect to the new 1-cells.

\subsection{Monoidal bicategories}

This section will show how the free monoidal bicategory on one object is monoidally biequivalent to the fundamental 2-groupoid of the coproduct
\[
\coprod_{k} B(k, \mathbb{R}^{1}).
\]
The proof is trivial as the free monoidal bicategory on one object is controlled by the coherence theorem for monoidal bicategories.

\begin{thm}[Coherence for monoidal bicategories I]
The strict functor
\[
\mathcal{F}X \rightarrow \mathcal{F}_{\textrm{Gr}}X
\]
induced by the universal property of $\mathcal{F}X$ is a monoidal biequivalence.
\end{thm}

This presentation of coherence for monoidal bicategories is just a special case of the coherence results for tricategories proven in \cite{g1}.  We now specialize to the case when $X$ is a singleton set.

\begin{thm}[Coherence for monoidal bicategories II]
The free monoidal bicategory on one object is monoidally biequivalent to the monoidal structure on
\[
\Pi_{2} \big( \coprod_{k} B(k, \mathbb{R}^{1}) \big)
\]
induced by the homotopy equivalences $B(k, \mathbb{R}^{1}) \simeq \mathcal{C}_{1}(k)/ \Sigma_{k}$.
\end{thm}
$\mb{Proof.}$
It is simple to compute that $\mathcal{C}_{1}(k)/ \Sigma_{k}$ is contractible, so that the bigroupoid $\Pi_{2} \big( \coprod_{k} B(k, \mathbb{R}^{1}) \big)$ has hom-categories which are all contractible.  The same holds for the free monoidal bicategory on one object by coherence \cite{gps,g1}, so we need only provide a monoidal functor $F$ that is biessentially surjective (that is, every object of the target is equivalent to one in the image) and which has the property that if $Fa \simeq Fb$ then $a \simeq b$.  Now the free monoidal bicategory on one object has a universal property, namely that strict functors from it to a given monoidal bicategory correspond to objects of that monoidal bicategory.  Thus to give a strict map from the free monoidal bicategory on one object to $\Pi_{2} \big( \coprod_{k} B(k, \mathbb{R}^{1}) \big)$, we need only choose an object of the target.

Consider the functor induced by the universal property sending the generating object $x$ of the free monoidal bicategory to the object $0 \in \mathbb{R} = B(1, \mathbb{R}^{1})$.  Since the monoidal structure on $\Pi_{2} \big( \coprod_{k} B(k, \mathbb{R}^{1}) \big)$ is induced from the operadic composition, the tensor product of $n$ copies of $0 \in \mathbb{R} = B(1, \mathbb{R}^{1})$ lands in the image of $\Pi_{2} B(n, \mathbb{R}^{1})$.  Since each space $B(n, \mathbb{R}^{1})$ is contractible, that means every object of $\Pi_{2} B(n, \mathbb{R}^{1})$ is equivalent to $0^{\otimes n}$, and thus the induced map $F$ from the free monoidal bicategory on one object is biessentially surjective.  We must also check that $Fa \simeq Fb$ implies $a \simeq b$, but this is also clear.  $\Box$

The proof above shows how coherence can provide alternate descriptions of free monoidal bicategories using topology.  The next section will begin the task of doing the converse:  using topological information to prove a coherence result.

\subsection{Braided monoidal bicategories}
Here we compare free braided monoidal bicategories with the fundamental 2-groupoid of the coproduct
\[
\coprod_{k} B(k, \mathbb{R}^{2}).
\]
Recall that for every algebra $X$ over the little 2-cubes operad, $\Pi_{2}X$ is a braided monoidal bicategory.  In particular, this is true of free algebras.  We have already seen that
\[
C_{2}(*) \simeq \coprod_{k} B(k, \mathbb{R}^{2}),
\]
where here * denotes a terminal space; we use this to fix the braided monoidal structure on $\Pi_{2} \big( \coprod_{k} B(k, \mathbb{R}^{2}) \big)$ by Theorem 15 and Corollary 16.  The main result is then the following, where $\mathcal{F}_{br}$ is the free braided monoidal bicategory functor.

\begin{thm}[Coherence for braided monoidal bicategories]\label{cohbmb}
The canonical map
\[
T:\mathcal{F}_{br}(*) \rightarrow \Pi_{2} \big( \coprod_{k} B(k, \mathbb{R}^{2}) \big)
\]
induced by sending the generating object $x$ to the point $(0,0) \in B(1, \mathbb{R}^{2}) = \mathbb{R}^{2}$ is a braided monoidal biequivalence.
\end{thm}

The proof of this theorem is based on results of Carter and Saito \cite{cs2} classifying surface braids in $\mathbb{R}^{4}$.  They show that two surface braids are equivalent if and only if they are related by a finite sequence of braid movie moves.  The braid movie moves give a completely algebraic description of the ambient isotopy relation, and we will see that the algebra they describe includes the braided monoidal bicategory axioms, for example see \cite{bl} for a construction of a braided monoidal 2-category of 2-tangles in 4-space using the braid movie moves.

Before beginning, we note that the coherence theorem for monoidal bicategories is used implicitly in this proof.  Thus we write, for example, tensors without bracketing as any two choices of brackets will be equivalent and any two equivalences between these choices will be uniquely isomorphic.

\textbf{Proof of \ref{cohbmb}.}
This is a braided monoidal functor since it is induced by the universal property, hence we must only show that it is a biequivalence.  Thus we must prove that $T$ is biessentially surjective and locally an equivalence of categories.

To show that $T$ is biessentially surjective, we must show that for every object $y$ in the target, there is an object $y'$ in the source for which $Ty'$ is equivalent to $y$.  Since all of the spaces $B(k, \mathbb{R}^{2})$ are connected, the equivalence class of an object in the target is determined completely by which space $B(k, \mathbb{R}^{2})$ contains the point $y$.  Thus to prove that the functor $T$ is biessentially surjective, we must show that, for every natural number $n$, there is some object $z_{n} \in \mathcal{F}_{br}*$ that maps to an object in the image of
\[
\Pi_{2} B(n, \mathbb{R}^{2}) \hookrightarrow \Pi_{2} \big( \coprod_{k} B(k, \mathbb{R}^{2}) \big)
\]
induced by the inclusion into the coproduct.  By definition, $T$ maps the generating object $x$ to $(0,0) \in B(1, \mathbb{R}^{2})$, so we have
\[
T(x^{\otimes n}) \simeq (Tx)^{\otimes n} = (0,0)^{\otimes n}.
\]
Now the monoidal structure in the target is that obtained by transfer from $\Pi_{2}\big( C_{2}(*) \big)$, and the monoidal structure there is given by the algebra structure of $C_{2}(*)$ over the operad $\mathcal{C}_{2}$.  By construction, the tensor product in $\Pi_{2}\big( C_{2}(*) \big)$ has the following property:  if $a_{i}$ is an object in $\Pi_{2} \mathcal{C}_{2}(n_{i})$ for $i=1, 2, \ldots, k$, then $a_{1} \otimes a_{2} \otimes \cdots \otimes a_{k}$ is an object in $\Pi_{2} \mathcal{C}_{2}(n_{1} + \cdots + n_{k})$.  This property is transferred back to $\Pi_{2} \big( \coprod_{k} B(k, \mathbb{R}^{2}) \big)$, so $(0,0)^{\otimes n}$ is an object in the essential image of $\Pi_{2} B(n, \mathbb{R}^{2})$.  Thus setting $z_{n} = x^{\otimes n}$ shows that $T$ is biessentially surjective.

Now we must show that $T$ induces an equivalence of categories on each hom-category, or in other words that $T$ is locally essentially surjective, locally full, and locally faithful.  Since every object in the free braided monoidal category on a single object $x$ is equivalent to $x^{k}$ for some value of $k$, we limit ourselves to this case.  Functoriality of $T$ will guarantee that $T$ is a local equivalence if it is true that $T$ is an equivalence on each hom-category of the form $\mathcal{F}_{br}(*)(x^{k}, x^{l})$.  This category is empty unless $k=l$, as is the category
\[
T(k,l):=\Pi_{2} \big( \coprod B(k, \mathbb{R}^{2}) \big)(T(x^{k}), T(x^{l})),
\]
so we are reduced to the case $k=l$.

To show that $T$ is locally essentially surjective when $k=l$, note that the objects in $T(k,k)$ are braids of $k$ strands, each of which is isomorphic to a composite of $\sigma_{i}$'s (the braid taking the $i$th strand over the $(i+1)$st).  By the braided monoidal functor axioms, $T$ sends the composite
\[
x^{k} \rightarrow x^{i-1}((xx)x^{k-i-1}) \stackrel{1(R1)}{\rightarrow}  x^{i-1}((xx)x^{k-i-1}) \rightarrow x^{k}
\]
to a braid isomorphic to $\sigma_{i}$, so $T$ is locally essentially surjective.

Finally, we must prove that $T$ is locally full and locally faithful.  Let $f,g \in \mathcal{F}_{br}(*)(x^k,x^k)$, and let $H$ be a 2-cell $Tf \Rightarrow Tg$, i.e.,  a homotopy between the maps $Tf, Tg$ fixing the boundaries, and so $Tf$ and $Tg$ represent equivalent braids in the braid group on $k$ letters.  Therefore $H$ is homotopic to a homotopy composed of a finite sequence of the braid group relations (see \cite{b})
\[
\begin{array}{cr}
\sigma_{i} \sigma_{j} = \sigma_{j} \sigma_{i}, & \quad |i - j| \geq 2 \\
\sigma_{i} \sigma_{i+1} \sigma_{i} = \sigma_{i+1} \sigma_{i} \sigma_{i+1} &
\end{array}
\]
transforming the braid $Tf$ into the braid $Tg$; here the homotopy corresponding to the first relation just reparametrizes the strings to alter the heights of the crossing points, while the homotopy corresponding to the second relation slides the crossing of strings $i+1$ and $i+2$ under the $i$th string.  Proving local fullness then reduces to showing that the homotopies corresponding to each of the relations above is in the image of $T$.

Now since $\otimes$ is a weak functor, we have a composite of functoriality isomorphisms
\[
(f \otimes 1) \circ (1 \otimes g) \cong f \otimes g \cong (1 \otimes g) \circ (f \otimes 1)
\]
for any pair of 1-cells $f,g$.  Consider the following 1-cells in $\mathcal{F}_{br}(*)(x^k,x^k)$.  Let $f$ be
\[
x \otimes x \otimes \cdots x \otimes x \otimes \cdots \otimes x \stackrel{1 \otimes \cdots \otimes 1 \otimes R_{x,x} \otimes 1 \otimes \cdots \otimes 1}{\longrightarrow} x \otimes x \otimes \cdots x \otimes x \otimes \cdots \otimes x
\]
where there are $k$ copies of the generator $x$, and $f$ consists of $k-2$ copies of the identity tensored together, along with one copy of $R_{x,x}$ which switches the $i$th copy of $x$ past the $(i+1)$st copy of $x$.  Define $g$ in the same way, but with the copy of $R_{x,x}$ switching the $j$th copy of $x$ past the $(j+1)$st.  Assume that $|i-j| \geq 2$.  Since $T$ is a strict functor, it sends the composite of functoriality isomorphisms above in $\mathcal{F}_{br}(*)$ to the composite of functoriality isomorphisms
\[
(Tf \otimes 1) \circ (1 \otimes Tg) \cong Tf \otimes Tg \cong (1 \otimes Tg) \circ (Tf \otimes 1).
\]
This functoriality isomorphism is in the homotopy class in $T(k,k)$ of the first braid group relation.  A similar proof shows that the second braid group relation is in the essential image using either of the 2-cells
\[
a \circ (R1) \circ a^{\centerdot} \circ (1R) \circ a \circ (R1) \Rightarrow (1R) \circ a \circ (R1) \circ a^{\centerdot} \circ (1R) \circ a
\]
appearing in the fourth braided monoidal bicategory axiom map; once again, this relies on the strictness of the induced map from the universal property.

Now we turn to showing that $T$ is locally faithful.  To prove this, assume that $H,K:f \Rightarrow g$ are mapped to the same homotopy class of maps in $T(k,k)$.  The homotopies $H,K$ are maps $D^{2} \rightarrow B_{e}(\mathbb{R}^{2}, k)$ which by Lemma \ref{transverse} (noting that these homotopies are necessarily disjoint from $\Sigma^{1}_{(k)}(\textrm{int }D^{2})$ by construction) give equivalent surface braids.  Thus there is a finite sequence of braid movie moves and locality changes which makes explicit the equivalence between these two surface braids.  If we can show that all of the braid movie moves and locality changes between such surface braids can be written in terms of the braided monoidal bicategory axioms, then $T$ is locally faithful since $TH = TK$ will imply that $H=K$.

Since we are only considering the graphs of homotopies between maps $I \rightarrow B(\mathbb{R}^{2}, k)$, the corresponding surface braids are covers which are not branched; thus we need only show that the C-I movie moves and locality changes can be written in terms of the braided monoidal bicategory axioms.  There are ten C-I moves, and we express each of them using the braided monoidal bicategory axioms individually.  In each case the braid movie move, written out here using the free monoid on the letters $s_{i}^{\pm 1}$, corresponds to an equivalence between between two surface braids, one presented to the left of the double-sided arrow and one presented to the right.  The task is then to show that the braided monoidal bicategory axioms imply that the 2-cell corresponding to the surface braid on the left is equal to the 2-cell corresponding to the surface braid on the right.  This is all straightforward, but we encourage the reader to consult Chapter 3 of \cite{cs2} for some useful illustrations.
\begin{itemize}
\item Move C-I-R1: This move is
\[
(s_{i}, s_{i}s_{j}s_{j}^{-1}, s_{j}s_{i}s_{j}^{-1}) \leftrightarrow (s_{i}, s_{j}s_{j}^{-1}s_{i}, s_{j}s_{i}s_{j}^{-1}),
\]
where $|i-j|>1$, and is a direct consequence of the functoriality of the tensor product.
\item Move C-I-R1': This move is
\[
(s_{i}, s_{i}s_{i}^{-1}s_{i}, s_{i}) \leftrightarrow (s_{i})
\]
and is one of the triangle identities for the braiding $R$ which is part of an adjoint equivalence.
\item Move C-I-R2: This move is
\[
(s_{i}s_{j}, s_{j}s_{i}, s_{i}s_{j}) \leftrightarrow (s_{i}s_{j}),
\]
where $|i-j|>1$, and is a direct consequence of the functoriality of the tensor product.
\item Move C-I-R3:  This move is
\[
(s_{i}s_{k}s_{j}, s_{k}s_{i}s_{j}, s_{k}s_{j}s_{i}, s_{j}s_{k}s_{i}) \leftrightarrow s_{i}s_{k}s_{j}, s_{i}s_{j}s_{k}, s_{j}s_{i}s_{k}, s_{j}s_{k}s_{i}),
\]
where  $|i-j|, |j-k|, |k-i|>1$, and follows from coherence for monoidal bicategories as the diagram consists only of coherence cells from the monoidal structure.
\item Move C-I-R4: This move is
\[
\begin{array}{l}
(s_{k}s_{i}s_{j}s_{i}, s_{i}s_{k}s_{j}s_{i}, s_{i}s_{j}s_{k}s_{i}, s_{i}s_{j}s_{i}s_{k}, s_{j}s_{i}s_{j}s_{k}) \leftrightarrow \\
\phantom{empty!!!} (s_{k}s_{i}s_{j}s_{i}, s_{k}s_{j}s_{i}s_{j}, s_{j}s_{k}s_{i}s_{j}s_{j}s_{i}s_{k}s_{j}, s_{j}s_{i}s_{j}s_{k}),
\end{array}
\]
where $|i-j|=1$ and $|i-k|, |j-k|>1$, and follows from functoriality of the tensor product.
\item Move C-I-M1: This move is
\[
(\textrm{empty word}) \leftrightarrow (\textrm{empty word}, s_{i}s_{i}^{-1}, \textrm{empty word})
\]
and follows from the invertibility of the 2-cell $1_{x \otimes x} \Rightarrow R_{x,x}^{\centerdot} \circ R_{x,x}$.
\item Move C-I-M2: This move is
\[
(s_{i}s_{i}^{-1}, \textrm{empty word}, s_{i}s_{i}^{-1}) \leftrightarrow (s_{i}s_{i}^{-1})
\]
and follows from the invertibility of the 2-cell $1_{x \otimes x} \Rightarrow R_{x,x}^{\centerdot} \circ R_{x,x}$.
\item Move C-I-M3:  This move is
\[
(s_{i}s_{j}s_{i}, s_{j}s_{i}s_{j}, s_{i}s_{j}s_{i}) \leftrightarrow (s_{i}s_{j}s_{i}),
\]
where $|i-j|=1$, and follows from the fourth braided monoidal bicategory axiom.
\item Move C-I-M4:  This move is
\[
\begin{array}{l}
(s_{i}s_{j}s_{k}s_{i}s_{j}s_{i}, s_{i}s_{j}s_{i}s_{k}s_{j}s_{i}, s_{j}s_{i}s_{j}s_{k}s_{j}s_{i}, s_{j}s_{i}s_{k}s_{j}s_{k}s_{i}, \\
\phantom{empty!!!} s_{j}s_{k}s_{i}s_{j}s_{k}s_{i}, s_{j}s_{k}s_{i}s_{j}s_{i}s_{k}, s_{j}s_{k}s_{j}s_{i}s_{j}s_{k}, s_{k}s_{j}s_{k}s_{i}s_{j}s_{k}) \leftrightarrow \\
(s_{i}s_{j}s_{k}s_{i}s_{j}s_{i}, s_{i}s_{j}s_{k}s_{i}s_{j}s_{i}, s_{i}s_{k}s_{j}s_{k}s_{i}s_{j}, s_{k}s_{i}s_{j}s_{k}s_{i}s_{j}, \\
\phantom{empty!!!} s_{k}s_{i}s_{j}s_{i}s_{k}s_{j}, s_{k}s_{j}s_{i}s_{j}s_{k}s_{j}, s_{k}s_{j}s_{i}s_{k}s_{k}s_{j}s_{k},  s_{k}s_{j}s_{k}s_{i}s_{j}s_{k}),
\end{array}
\]
where $k=j+1 = i+2$ or $k= j-1 = i-2$, and follows from the third and fourth axioms, together with the fact that $R_{(-|-,-)}, R_{(-,-|-)}$ are modifications and $R$ is a pseudonatural transformation.
\item Move C-I-M5: This move is
\[
(s_{j}s_{i}, s_{i}^{-1}s_{i}s_{j}s_{i}, s_{i}^{-1}s_{j}s_{i}s_{j}) \leftrightarrow (s_{j}s_{i}, s_{j}s_{i}s_{j}^{-1}s_{j}, s_{i}^{-1}s_{j}s_{i}s_{j}),
\]
where $|i-j|=1$, and follows from the fourth braided monoidal bicategory axiom together with the fact that the mate of the naturality square for $R$
\[
\xy
{\ar^{R_{y,x}1} (0,0)*+{yxz}; (0,15)*+{xyz} };
{\ar^{R_{xy,z}} (0,15)*+{xyz}; (25,15)*+{zxy} };
{\ar_{R_{yx,z}} (0,0)*+{yxz}; (25,0)*+{zyx} };
{\ar_{1R_{y,x}} (25,0)*+{zyx}; (25,15)*+{zxy} };
(12.5, 7.5)*{\cong}
\endxy
\]
under the adjoint equivalence $R \dashv R^{\centerdot}$ is the naturality square for $R$ shown below.
\[
\xy
{\ar_{R_{y,x}^{\centerdot}1} (0,15)*+{xyz}; (0,0)*+{yxz} };
{\ar^{R_{xy,z}} (0,15)*+{xyz}; (25,15)*+{zxy} };
{\ar_{R_{yx,z}} (0,0)*+{yxz}; (25,0)*+{zyx} };
{\ar^{1R_{y,x}^{\centerdot}} (25,15)*+{zxy}; (25,0)*+{zyx} };
(12.5, 7.5)*{\cong}
\endxy
\]
\item Locality changes:  First, recall that elementary braid changes between surface braids are 2-cells in the braided monoidal bicategory $\Pi_{2} \Big( \coprod_{k} B(k, \mathbb{R}^{2}) \Big)$.  Second, locality changes involve writing some terms of a braid movie as $w_{i} = u_{i}v_{i}$ which corresponds to composing 1-cells in the braided monoidal bicategory.  Thus the claim that two braid movies which differ by a locality change are equivalent follows from the naturality of the \textbf{Gray}-category structure isomorphisms with respect to 2-cells. $\Box$
\end{itemize}

\begin{cor}\label{uniqueiso}
In $\mathcal{F}_{br}(*)$, there is at most one isomorphism between any pair of parallel 1-cells.
\end{cor}
$\mb{Proof.}$
Since $\pi_{2}B(k,\mathbb{R}^{2}) = 0$ for all $k$, there is at most a single isomorphism $1_{a} \Rightarrow 1_{a}$ for any object $a$ in $\mathcal{F}_{br}(*)$.  But since $\mathcal{F}_{br}(*)$ is a bigroupoid, this implies that the same holds for any parallel pair of 1-cells.  $\Box$

In the next section, we require versions of Theorem \ref{cohbmb} and Corollary \ref{uniqueiso} for the free braided monoidal bicategory on many objects.  The proofs of these two results are the same as the previous results with the addition that all of the homotopies are additionally labeled by a set $S$.  Since the labels do not affect the geometry at all, requiring a geometric condition to hold (such as two surface braids being equivalent) in the presence of labels is logically equivalent to it holding without labels and requiring that all constructions preserve labels.  For example, a path in the labeled configuration space from $(x_{i},s_{i})$ to $(y_{i},t_{i})$ is the same as
\begin{itemize}
\item the statement that $s_{i}=t_{i}$ if the path connects $x_{i}$ to $y_{i}$,
\item and a path $\gamma$ in the unlabeled configuration space.
\end{itemize}
Thus all of the geometric and topological results required can be immediately generalized to the case of a set of labels.

\begin{thm}\label{cohbmb2}
Let $S$ be a set, seen as a discrete bicategory.  The free braided monoidal bicategory on $S$ is braided monoidally biequivalent to $\Pi_{2} \Big( C_{2}(S) \Big)$, and the map
\[
T: \mathcal{F}_{br}(S) \rightarrow \Pi_{2} \Big( \coprod_{k} B(k, \mathbb{R}^{2}; S) \Big)
\]
induced by the universal property sending the element $s \in S$ to the point $\Big( (0,0); s \Big) \in B(1, \mathbb{R}^{2}; S)$ is a braided monoidal biequivalence, where the target is given the lifted structure.
\end{thm}

\begin{cor}\label{uniqueiso2}
In $\mathcal{F}_{br}(S)$, there is at most one isomorphism between any pair of parallel 1-cells.
\end{cor}
$\mb{Proof.}$  We once again show that $\pi_{2}$ of some space is zero, in this case that space is $B(k, \mathbb{R}^{2}; S)$.  We have the fibration
$B(k, \mathbb{R}^{2}; S) \rightarrow B(k, \mathbb{R}^{2})$ with fiber the discrete space $S^{k}$, so by the long exact sequence in homotopy groups this map induces an isomorphism on $\pi_{i}$ for $i>1$, and so the results follows from Corollary \ref{uniqueiso}.  $\Box$

\section{Strictification}
This section presents an alternate coherence theorem for braided monoidal bicategories that we prove using the results of Section 4.3.  It is a standard strictification theorem and states that every braided monoidal bicategory is braided monoidally biequivalent to a braided monoidal 2-category in the sense of Crans \cite{cr}.  This coherence theorem relates the work here to previous work on braided monoidal 2-categories by Kapranov and Voevodsky \cite{kv2}, Baez and Neuchl \cite{bn}, Baez and Langford \cite{bl}, Crans \cite{cr}, and Day and Street \cite{ds}.

The key ingredients for the proof of this theorem are Corollary \ref{uniqueiso2} and the strictification theorem for monoidal bicategories.  We shall assume that the reader is familiar with this strictification $B \mapsto \textrm{Gr}B$ as presented in \cite{g1}, specialized to the case when the tricategory in question only has one object and is thus a monoidal bicategory.  Roughly speaking, $\textrm{Gr}B$ has objects which are strings of objects in $B$, 1-cells which are strings of 1-cells in $B$, and 2-cells between strings which are 2-cells between their composites in $B$.

It should be noted that the result in this section is not an immediate consequence of the existence of a strictification, but rather a consequence of the construction of the strictification $\textrm{Gr} B$.  Coherence for tricategories in its one-object form implies that any braided monoidal bicategory is monoidally biequivalent to a $\mb{Gray}$-monoid equipped with a braided structure; this is merely using coherence to strictify the underlying monoidal bicategory and then lifting the braided structure.  The important thing to remember is that a braided monoidal 2-category is not just a $\mb{Gray}$-monoid with a braiding, but is also required to satisfy the additional unit axioms added by Crans.  The content of this section is then that it is always possible to equip the particular strictification $\textrm{Gr} B$ with the full structure of a braided monoidal 2-category.

\begin{thm}[Strictification for braided monoidal bicategories]\label{strictbmb}
Let $B$ be a braided monoidal bicategory.  Then the strictification $\textrm{Gr} B$ can be equipped with the structure of a braided monoidal 2-category, and the monoidal biequivalences
\[
\begin{array}{c}
f:B \rightarrow \textrm{Gr}B \\
e:\textrm{Gr}B \rightarrow B
\end{array}
\]
can each be equipped with the structure of a braided monoidal functor.
\end{thm}
$\mb{Proof.}$
First we give $\textrm{Gr}B$ the structure of a braided monoidal 2-category.  It is already a $\mb{Gray}$-monoid by construction, so we need only define the braided structure making sure that it satisfies the conditions listed above.

Recall that an object of $\textrm{Gr}B$ is a string $X = (X_{n}, X_{n-1}, \ldots, X_{1})$ of objects of $B$; we also include the empty string as an object.  For any such string, we have the object
\[
e(X) = ( \cdots (X_{n} \otimes X_{n-1}) \otimes X_{n-2}) \otimes \cdots \otimes X_{2}) \otimes X_{1}
\]
in $B$.  Given two non-empty strings $X,Y$, we define $R_{X,Y}$ to be the 1-cell given by
\[
R_{e(X), e(Y)}: e(X) \otimes e(Y) \rightarrow e(Y) \otimes e(X)
\]
in $B$.  If either $X$ or $Y$ is the empty string, we define $R_{X,Y}$ to be the identity 1-cell in $\textrm{Gr}B$.  A similar definition is made for $R_{X,Y}^{\centerdot}$; the unit and counit of this adjoint equivalence are then inherited from the unit and counit in $B$.  It is then clear that the first two unit conditions are satisfied.

Now we must define isomorphisms $R_{(X|Y,Z)}$ and $R_{(X,Y|Z)}$.  By Corollary \ref{uniqueiso} and the construction of $\textrm{Gr}B$, both of these isomorphisms are uniquely determined by the braided monoidal structure on $B$ as follows.  Since $\textrm{Gr}B$ is a monoidal 2-category, the 2-cell $R_{(X|Y,Z)}$ has source and target
\[
(1_{Y} \otimes R_{X,Z}) \circ (R_{X,Y} \otimes 1_{Z}) \Rightarrow R_{X, YZ}.
\]
Now a 2-cell $\alpha$ in $\textrm{Gr}B$ is just a 2-cell in $B$ with source given by applying the functor $e$ to the source of $\alpha$ and target given by applying $e$ to the target of $\alpha$.  The functor $e$ applied to the 1-cell $(1_{Y} \otimes R_{X,Z}) \circ (R_{X,Y} \otimes 1_{Z})$ in $\textrm{Gr}B$ consists of a composite of 1-cells in $B$ all of which arise from the braided monoidal structure; the same holds when applying $e$ to $R_{X, YZ}$.  In particular, this shows that the 1-cell source and target of the 2-cell $R_{(X|Y,Z)}$ in $\textrm{Gr}B$ are both cells which arise from the free braided monoidal bicategory on objects $X,Y,Z$, so there is a unique 2-cell coherence isomorphism between them in $B$.  This unique cell is the 2-cell $R_{(X|Y,Z)}$ in $\textrm{Gr}B$.  Additionally, $\textrm{Gr} B$ satisfies the axioms for a braided monoidal bicategory by the same reasoning as in each case the two different pastings which must be shown to be equal for a given axiom to hold are both 2-cells in a free braided monoidal bicategory with the same source and target, hence are equal.

It only remains to check the last six unit conditions to show that $\textrm{Gr} B$ is a braided monoidal 2-category.  Thus we must show that some of these uniquely determined isomorphisms are in fact the identity.  In each case, we need only show that the 2-cell in question has the same source as target; if this is true, then the identity is a valid candidate for the coherence 2-cell, and so by uniqueness must be.  This is trivial by examining which instances of $R$ are actually identities and using the fact that $\textrm{Gr} B$ is a $\mb{Gray}$-monoid.  As an example, consider $R_{(A|B,I)}$.  The 1-cell source of this 2-cell is $R_{A,BI}$ which is equal to $R_{A,B}$ since $I$ is a strict unit.  The 1-cell target of this 2-cell is $(1_{B} \otimes R_{A,I}) \circ (R_{A,B} \otimes 1_{I})$.  By definition, $R_{A,I}$ is the identity on $A$, so $1_{B} \otimes R_{A,I}$ is $1_{BA}$.  Additionally, $R_{A,B} \otimes 1_{I} = R_{A,B}$, so the target becomes
\[
(1_{B} \otimes R_{A,I}) \circ (R_{A,B} \otimes 1_{I}) = 1_{AB} \circ R_{A,B} = R_{A,B},
\]
and therefore $R_{(A|B,I)}$ must be the identity 2-cell by coherence.

It is easy to show that $e$ and $f$ can be extended to braided monoidal functors.  These are already monoidal functors, so the only thing left to define is the invertible modification $U$ and then check two axioms.  For the functor $e$, $U$ has its component at $X,Y$ a 2-cell of the shape shown below.
\[
\xy
{\ar^{R} (0,0)*+{e(X)e(Y)}; (30,0)*+{e(Y)e(X)} };
{\ar^{\chi} (30,0)*+{e(Y)e(X)}; (30,-12)*+{e(YX)} };
{\ar_{\chi} (0,0)*+{e(X)e(Y)}; (0,-12)*+{e(XY)} };
{\ar_{e(R)} (0,-12)*+{e(XY)}; (30,-12)*+{e(YX)} };
(15,-6)*{\Downarrow U}
\endxy
\]
The 1-cells $\chi$ are associativity (if neither $X$ nor $Y$ is the unit) or unit (if at least one of $X$ or $Y$ is the unit) constraints from the monoidal structure, so we define $U$ to be the unique coherence 2-cell by Corollary \ref{uniqueiso2}.  The two axioms then follow immediately by coherence.  The construction for $f$ is analogous. $\Box$

\begin{cor}
Every braided monoidal 2-category in the sense of Baez-Neuchl is braided monoidally biequivalent to a braided monoidal 2-category in the sense of Crans.
\end{cor}

\providecommand{\bysame}{\leavevmode\hbox to3em{\hrulefill}\thinspace}
\providecommand{\MR}{\relax\ifhmode\unskip\space\fi MR }
\providecommand{\MRhref}[2]{%
  \href{http://www.ams.org/mathscinet-getitem?mr=#1}{#2}
}
\providecommand{\href}[2]{#2}

\end{document}